\documentclass[a4paper]{amsart}

\usepackage{amssymb}
\usepackage{amscd}

\usepackage{times}

\usepackage{t1enc}

\usepackage[english,francais]{babel}

\def\makeatother{\catcode64=\active}

\newcommand{\UseRsfsAsCalli}{
\DeclareFontFamily{U}{rsfs}{}
\DeclareFontShape{U}{rsfs}{m}{n}{%
   <5>rsfs5%
   <6>rsfs10%
   <7>rsfs7%
   <8>rsfs10%
   <9>rsfs10%
   <10>rsfs10%
   <11>rsfs10%
   <12>rsfs10%
   <14>rsfs10%
   <17>rsfs10%
   <20>rsfs10%
   <25>rsfs10}{}
\DeclareMathAlphabet{\callig}{U}{rsfs}{m}{n}
\newcommand{\calli}[1]{{\callig ##1\/}}
}

\UseRsfsAsCalli

\DeclareFontFamily{U}{wncy}{}
\DeclareFontShape{U}{wncy}{m}{n}{%
   <5>wncyr5%
   <6>wncyr6%
   <7>wncyr7%
   <8>wncyr8%
   <9>wncyr9%
   <10>wncyr10%
   <11>wncyr10%
   <12>wncyr6%
   <14>wncyr7%
   <17>wncyr8%
   <20>wncyr10%
   <25>wncyr10}{}
\DeclareMathAlphabet{\cyrille}{U}{wncy}{m}{n}

\DeclareMathAlphabet{\boldita}{U}{cmr}{bx}{it}
\newcommand{\boldit}[1]{{\boldita #1}}

\DeclareMathAlphabet{\eulercal}{U}{eus}{m}{n}
\DeclareMathAlphabet{\eulerrm}{U}{eur}{m}{n}

\newbox\dummybox
\def\mysubscripts{
\setbox\dummybox\hbox{$$\fontdimen16\textfont2=2.8pt$$
$$\fontdimen16\scriptfont2=1.9pt$$}}
\mysubscripts

\def\greeksubscript#1{{#1}}

\newskip\placeabove
\newcommand{\noqed}{\renewcommand{\qed}{}}

\newcommand{\numero}{$\hbox{n}^{\hbox{\footnotesize o}}$ }

\newcommand{\varleq}{\leqslant}
\newcommand{\vargeq}{\geqslant}

\catcode`@=11
\def\frownfill{$\m@th
\braceld\leaders\vrule\hfill\bracerd$}
\def\overfrown#1{\mathop{\vbox{\ialign{##\crcr\noalign{\kern3pt}
\frownfill\crcr\noalign{\kern3pt\nointerlineskip}
$\hfil\displaystyle{#1}\hfil$\crcr}}}\limits}
\def\interieur#1{\mathord{\vbox{\ialign{##\crcr\noalign{\kern3pt}
$\hfill\circ\hfill$\crcr\noalign{\kern1pt\nointerlineskip}
\frownfill\crcr\noalign{\kern3pt\nointerlineskip}
$\hfil\displaystyle{#1}\hfil$\crcr}}}}
\makeatother

\newbox\oublieux
\newcommand{\forgetheight}[1]{\setbox\oublieux=\hbox{#1}%
\ht\oublieux=0pt%
\box\oublieux
}

\numberwithin{equation}{subsection}

\newcommand{\clap}[1]{\hbox to 0pt{\hss #1\hss}}
\newcommand{\mathclap}[1]{{\mathchoice{\clap{$\displaystyle #1$}}
{\clap{$#1$}}
{\clap{$\scriptstyle #1$}}
{\clap{$\scriptscriptstyle #1$}}}}

\renewcommand{\atop}[2]{\genfrac{}{}{0pt}{}{#1}{#2}}

\newbox\isobox
\newdimen\isodim
\newdimen\isoscriptdim
\newdimen\isoscriptscriptdim
\newcommand{\iso}{\mathrel{
\setbox\isobox\hbox{$\longrightarrow$}
\isodim=\wd\isobox
\setbox\isobox\hbox{$\scriptstyle\longrightarrow$}
\isoscriptdim=\wd\isobox
\setbox\isobox\hbox{$\scriptscriptstyle\longrightarrow$}
\isoscriptscriptdim=\wd\isobox
\mathchoice
{\hbox to\isodim{\hfil\lower 0.2ex\clap{$\widetilde{}$}%
\clap{$\longrightarrow$}\hfil}}%
{\hbox to\isodim{\hfil\lower 0.2ex\clap{$\widetilde{}$}%
\clap{$\longrightarrow$}\hfil}}%
{\hbox to\isoscriptdim{\hfil\lower 0.50ex%
\clap{$\scriptstyle\widetilde{}$}%
\clap{$\scriptstyle\longrightarrow$}\hfil}}%
{\hbox to\isoscriptscriptdim{\hfil\lower 0.60ex%
\clap{$\tilde{}$}%
\clap{$\scriptscriptstyle\longrightarrow$}\hfil}}}}

\catcode`@=11
\atdef@ surj{\ampersand@
  \ifCD@ \global\bigaw@\minCDarrowwidth \else \global\bigaw@\minaw@ \fi
 \ifCD@\enskip\fi
   \mathrel{\mathop{\hbox to\bigaw@{\rightarrowfill@\displaystyle
     \kern-2ex{$\to$}}}}%
 \ifCD@\enskip\fi \ampersand@}
\atdef@ i#1i#2i{\ampersand@
  \ifCD@ \global\bigaw@\minCDarrowwidth \else 
  \global\bigaw@\minaw@\fi
  \setbox\isobox\hbox{$\widetilde{\hbox{}}$}%
  \ifdim\wd\isobox>\bigaw@\global\bigaw@\wd\isobox\fi
  \setboxz@h{$\m@th\scriptstyle\;{#1}\;\;$}%
  \ifdim\wdz@>\bigaw@\global\bigaw@\wdz@\fi
  \@ifnotempty{#2}{\setbox\@ne\hbox{$\m@th\scriptstyle\;{#2}\;\;$}%
    \ifdim\wd\@ne>\bigaw@\global\bigaw@\wd\@ne\fi}%
  \ifCD@\enskip\fi
  \isodim\bigaw@
  \divide\isodim by 2
  \mathrel{\mathop{\hbox to\bigaw@{\hfil
        \clap{\hbox to \bigaw@{\rightarrowfill@\displaystyle}}%
        \lower 0.5ex\clap{$\widetilde{\hbox to\isodim{\hfil}}$}%
        \hfil}}%
    \limits^{#1}\@ifnotempty{#2}{_{#2}}}%
  \ifCD@\enskip\fi \ampersand@}
\makeatother

\CDat

\catcode`@=11
\atdef@ @{@}
\makeatother

\newcommand{\inverselimit}{\mathop{\underleftarrow%
{\vbox{\hrule width 0pt height 0pt depth 2pt}\lim}}}

\newcommand{\mut}[2]{\if #1\infty
\if #20 \mathbf Q/\mathbf Z
\else \mathbf Q/\mathbf Z(#2)
\fi\else
\if #20 \mathbf Z/#1\mathbf Z \else\if #21 \mu_{#1} \else
\mu^{\otimes#2}_{#1}\fi\fi\fi}

\newcommand{\textinmath}[1]{{\mathchoice%
{\hbox{\fontshape{n}\selectfont #1}}%
{\hbox{\fontshape{n}\selectfont #1}}%
{\hbox{\fontshape{n}\selectfont\scriptsize #1}}%
{\hbox{\fontshape{n}\selectfont\tiny #1}}}}

\DeclareMathOperator{\Det}{Det}
\DeclareMathOperator{\Tr}{Tr}
\renewcommand{\ker}{\mathop{\textinmath{Ker}}}

\DeclareMathOperator{\Br}{Br}
\DeclareMathOperator{\Pic}{Pic}

\DeclareMathOperator{\NS}{NS}

\DeclareMathOperator{\Supp}{Supp}
\DeclareMathOperator{\Fr}{Fr}
\DeclareMathOperator{\Gal}{Gal}
\DeclareMathOperator{\Spec}{Spec}

\DeclareMathOperator{\Hom}{Hom}

\DeclareMathOperator{\Lie}{Lie}

\DeclareMathOperator{\rg}{rg}

\DeclareMathOperator{\Jac}{Jac}

\DeclareMathOperator{\reel}{Re}
\DeclareMathOperator{\Res}{Res}

\newcommand{\eff}{_\textinmath{eff}}
\newcommand{\etale}{_\textinmath{\'et}}

\newcounter{lesrems}
\renewcommand{\thelesrems}{\roman{lesrems}}
\newenvironment{listrems}{\bgroup
\setcounter{lesrems}{0}
\renewcommand{\item}{{\refstepcounter{lesrems}%
\ifnum\value{lesrems}>1\par\fi
\fontshape{n}\selectfont (\thelesrems) }}%
}{\egroup}

\newenvironment{enumerer}{\bgroup
\setcounter{lesrems}{0}%
\begin{list}{(\roman{lesrems})}{%
\labelwidth=3em
\itemsep=0pt
\usecounter{lesrems}}}
{\end{list}\egroup}

\newtheorem{theo}{Th\'eor\`eme}[subsection]
\newtheorem{lemme}[theo]{Lemme}
\newtheorem{prop}[theo]{Proposition}
\newtheorem{cor}[theo]{Corollaire}

\theoremstyle{definition}
\newtheorem{defi}{D\'efinition}[subsection]
\newtheorem{defis}[defi]{D\'efinitions}

\theoremstyle{remark}
\newtheorem{rem}[theo]{Remarque}
\newtheorem{rems}[theo]{Remarques}

\newtheorem{nota}[defi]{Notation}
\newtheorem{notas}[defi]{Notations}
\newtheorem{conv}[defi]{Convention}
\newtheorem{hypo}{Hypoth\`ese}

\newtheorem{hyposgeo}[hypo]{Hypoth\`eses g\'eom\'etriques}
\newtheorem{hypoari}[hypo]{Hypoth\`ese arithm\'etique}

\catcode`@=11
\def\merci{\normalfont\small  \skip@28\p@ \advance\skip@-\lastskip
  \advance\skip@-\baselineskip \vskip\skip@
  \vtop \bgroup
}
\def\endmerci{
  \egroup
  \skip@32\p@\@plus 14\p@ \advance\skip@-\baselineskip
  \vskip\skip@}
\makeatother

\def\jourdhui{\number\day\space \ifcase\month\or
Janvier\or F\'evrier\or Mars\or
Avril\or Mai\or Juin\or Juillet\or Ao\^ut\or Septembre\or
Octobre\or Novembre\or D\'ecembre\fi
\space\number\year}

\catcode`@=11
\def\add@accent#1#2{{%
   \setbox\@tempboxa\hbox{#2}%
   \accent#1 #2}}
\makeatother

\newcommand{\cardinal}{\sharp}

\newcommand{\ZZ}{{\mathbf Z}}
\newcommand{\QQ}{{\mathbf Q}}
\newcommand{\RR}{{\mathbf R}}
\newcommand{\Rplus}{{\RR_{\geq 0}}}
\newcommand{\Rplusetoile}{{\RR_{> 0}}}
\newcommand{\unreel}{B}
\newcommand{\CC}{{\mathbf C}}

\newcommand{\Ll}{\ell}

\newcommand{\Cont}{\calli C^0}

\newcommand{\sous}{\backslash}

\newcommand{\etoile}{^*}
\newcommand{\ff}{F}
\newcommand{\ee}{E}
\newcommand{\normeEaF}{N_{\ee/\ff}}
\newcommand{\cloture}[1]{\overline{#1}}
\newcommand{\fbarre}{{\cloture{\ff}}}
\newcommand{\sep}{^s}
\newcommand{\fs}{{\ff\sep}}

\newcommand{\FF}{{\mathbf F}}

\newcommand{\qq}{q}
\newcommand{\car}{p}	
\newcommand{\qf}{q_\ff}
\newcommand{\Fqf}{\FF_{\qf}}
\newcommand{\anneau}{{\calli O}}
\newcommand{\hs}{^{\text{hs}}}
\newcommand{\SSS}{S}
\newcommand{\cc}{{\calli C}}

\newcommand{\genre}{g}
\newcommand{\gf}{\genre_\ff}
\newcommand{\hf}{h_\ff}
\newcommand{\adeles}{{\boldit A_\ff}}

\newcommand{\placesde}[1]{{M_{#1}}}
\newcommand{\Mf}{{\placesde\ff}}

\newcommand{\pp}{{\mathfrak p}}
\newcommand{\pP}{{\mathfrak P}}
\newcommand{\fp}{{\ff_\pp}}
\newcommand{\fpbarre}{{\overline\ff_\pp}}
\newcommand{\Op}{\anneau_\pp}
\newcommand{\OP}{\anneau_\pP}
\newcommand{\Fp}{\FF_\pp}

\newcommand{\Fpbarre}{{\overline\FF_\pp}}
\newcommand{\FPbarre}{{\overline\FF_\pP}}
\newcommand{\normede}[1]{\vert#1\vert_\pp}
\newcommand{\norme}{\normede\cdot}
\newcommand{\vp}{{v_\pp}}
\newcommand{\Haar}[1]{{\mathrm d}#1\,}
\newcommand{\dxp}{\Haar {x_\pp}}

\newcommand{\VV}{V}
\newcommand{\Vbarre}{{\overline\VV}}
\newcommand{\Vs}{{\VV^s}}
\newcommand{\UU}{U}                   
\newcommand{\ferme}{K}
\newcommand{\modelede}[1]{{\calli #1}}
\newcommand{\mV}{\modelede \VV}
\newcommand{\structural}{{\calli O}}
\newcommand{\OV}{\structural_{\VV}}
\newcommand{\OVx}{\structural_{\VV,x}}
\newcommand{\canoniquede}[1]{\omega_{#1}}
\newcommand{\omegaV}{\canoniquede\VV}
\newcommand{\antican}{\canoniquede\VV^{-1}}
\newcommand{\Ceffde}[1]{C\eff(#1)}
\newcommand{\CeffV}{\Ceffde\VV}
\newcommand{\CeffVbarre}{\Ceffde\Vbarre}
\newcommand{\LL}{L} 
\newcommand{\mL}{\modelede\LL}
\newcommand{\sectiondeL}{{\boldit s}}
\newcommand{\dual}{^\vee}

\newcommand{\PP}{{\mathbf P}}         
\newcommand{\affine}{{\mathbf A}}

\newcommand{\metriquede}[1]{\Vert#1\Vert_\pp}
\newcommand{\metrique}{\metriquede\cdot}
\newcommand{\metriques}{(\metrique)_{\pp\in\Mf}}

\newcommand{\hauteur}{H}
\newcommand{\hauteurs}{{\calli H}}
\newcommand{\oubli}{{\boldit o}}
\newcommand{\HH}{{\boldit H}}
\newcommand{\HK}{{\boldit H}_K}
\newcommand{\HPp}{H_{P,\mathfrak p}}
\newcommand{\HP}{H_{P}}
\newcommand{\HPZ}{H_{P_0}}
\newcommand{\nHde}[1]{{n_{#1,H}}}
\newcommand{\nWH}{{\nHde W}}
\newcommand{\nUH}{{\nHde U}}
\newcommand{\aLde}[1]{{a_{#1}(\LL)}}
\newcommand{\zetaHde}[1]{{\zeta_{#1,\HH}}}
\newcommand{\zetaHKde}[1]{{\zeta_{#1,\HK}}}
\newcommand{\zetaHU}{{\zetaHde U}}
\newcommand{\zetaHV}{{\zetaHde V}}
\newcommand{\zetaHKV}{{\zetaHKde V}}

\newcommand{\mesure}{\boldsymbol\omega}
\newcommand{\mesp}{{\mesure_\pp}}
\newcommand{\coeffconv}{\lambda_\pp}
\newcommand{\mesH}{\mesure_\hauteur}
\newcommand{\tauV}{\tau_\hauteur(\VV)}

\newcommand{\alphaV}{\alpha(\VV)}
\newcommand{\alphaetoileV}{\alpha^*(\VV)}
\newcommand{\betaV}{\beta(\VV)}
\newcommand{\thetaV}{\theta_\hauteur(\VV)}
\newcommand{\thetaetoileV}{\theta^*_\hauteur(\VV)}
\newcommand{\caracteristiquede}[1]{\chi_\greeksubscript{#1}}
\newcommand{\chiCeff}{{\caracteristiquede\CeffV}}

\newcommand{\compact}{K}
\newcommand{\rad}{{\calli R}}
\newcommand{\raduni}{\rad_u}
\newcommand{\norma}{\calli N}
\newcommand{\caract}{X^*}
\newcommand{\cocaract}{X_*}
\newcommand{\caractC}{{\mathfrak a}}
\newcommand{\caractCZ}{{\mathfrak a}_0}
\newcommand{\Gmde}[1]{\mathbf G_{m,#1}}
\newcommand{\relatifa}[2]{\,\vphantom{#2}_{#1}#2}
\newcommand{\relatif}{\relatifa\ff}
\newcommand{\relatifp}{\relatifa\fp}
\newcommand{\racines}{\Phi}
\newcommand{\Fracines}{\relatif\racines}
\newcommand{\Fpracines}{\relatifp\racines}
\newcommand{\base}{\Delta}
\newcommand{\Fbase}{\relatif\base}
\newcommand{\Fpbase}{\relatifp\base}
\newcommand{\positive}{\racines^+}
\newcommand{\Fpositive}{\relatif\positive}
\newcommand{\Weyl}{W}
\newcommand{\FWeyl}{\relatif\Weyl}
\newcommand{\corac}{\check}
\newcommand{\poids}[1]{\varpi_{#1}}
\newcommand{\frakr}{\mathfrak r}
\newcommand{\demiracines}[1]{\rho_\greeksubscript{#1}}

\newcommand{\Eisen}[2]{{E^{#1}_{#2}}}
\newcommand{\EGP}{{\Eisen GP}}
\newcommand{\EPPZ}{{\Eisen P{P_0}}}
\newcommand{\EGPZ}{{\Eisen G{P_0}}}

\def\forget#1{}

\title[Points de hauteur born\forgetheight{\'e}e en caract\forgetheight{\'e}ristique finie]
{Points de hauteur born\'ee\\
sur les vari\'et\'es de drapeaux\\
en caract\'eristique finie}
\author{Emmanuel Peyre}
\address{Institut Fourier\\
UFR de Math\'ematiques, UMR 5582\\
Universit\'e de Grenoble I et CNRS\\
BP 74\\ 38402 Saint-Martin d'H\`eres CEDEX\\ France}
\urladdr{http://www-fourier.ujf-grenoble.fr/\~{}peyre}
\email{Emmanuel.Peyre@@ujf-grenoble.fr}

\begin{document}

\begin{abstract}
Le but de cet article est d'appliquer les travaux de Morris 
sur les s\'eries d'Eisentein en caract\'eristique finie
\`a l'\'etude asymptotique des points rationnels
de hauteur born\'ee sur une vari\'et\'e de drapeaux
g\'en\'e\-ra\-li\-s\'ee obtenue comme quotient d'un groupe
alg\'ebrique semi-simple au-dessus d'un corps global
de caract\'eristique finie par un sous-groupe para\-bo\-li\-que r\'eduit.
La formule obtenue pour le r\'esidu de la fonction z\^eta des hauteurs
a une interpr\'etation analogue \`a celle connue pour un corps de nombres.
\end{abstract}
\forget{\begin{altabstract}
The aim of this paper is to apply the work of Morris
on Eisenstein series over global function fields
to the study of the asymptotic behavior of the points of
bounded height on a generalized flag variety defined as
the quotient of a semi-simple algebraic group
by a reduced parabolic subgroup over such a field.
The formula obtained for the height zeta function has an interpretation
similar to the one known over a number field.
\end{altabstract}}
\maketitle

\tableofcontents
\section*{Introduction}
La compr\'ehension du comportement asymptotique des points rationnels
de hauteur born\'ee sur les vari\'et\'es presque de Fano au-dessus
d'un corps de nombres a fortement progress\'e ces derni\`eres ann\'ees
notamment gr\^ace \`a l'impulsion donn\'ee par Manin
(cf. \cite{batyrevmanin:hauteur}, \cite{fmt:fano},
\cite{peyre:fano}, \cite{salberger:tamagawa} et 
\cite{batyrevtschinkel:tamagawa}). Il serait naturel que
le formalisme d\'evelopp\'e dans ce cadre s'\'etende dans
une certaine mesure au cas des corps globaux de
caract\'eristique finie. Il \'etait donc tentant de chercher
une formule pour le r\'esidu de la fonction z\^eta des hauteurs
pour les vari\'et\'es de drapeaux g\'en\'eralis\'ees sur un tel corps. Deux
raisons motivent cet exemple; tout d'abord de telles
formules asymptotiques ont \'et\'e obtenues dans des cas particuliers
(cf. \cite[\S2.5 in fine]{serre:mordellweil}, \cite{hsia:dynamical}),
d'autre part le r\^ole jou\'e par les travaux de Langlands 
dans la d\'emonstration de ces formules asymptotiques 
pour les vari\'et\'es de drapeaux sur un corps de nombres
(\cite{fmt:fano}, \cite{peyre:fano}) peut \^etre jou\'e par 
ceux de Morris dans le cas d'un corps de fonctions global 
(cf. \cite{morris:eisenstein:cusp}  et \cite{morris:eisenstein:general}).
\par
Entre la premi\`ere version de ce texte et sa soumission,
d'autres auteurs ont fait progresser cette extension au cadre fonctionnel
du programme initi\'e par Manin. D'une part, King Fai Lai et
Kit Ming Yeung dans \cite{laiyeung:flag}, \'ecrit
ind\'ependamment de notre texte, se sont \'egalement
int\'eress\'es aux vari\'et\'es de drapeaux dans le cadre fonctionnel,
sans toutefois donner d'interpr\'etation de la constante, ce qui constitue
le point crucial de notre travail. D'autre part D. Bourqui, a trait\'e 
de mani\`ere compl\`ete le cas d\'elicat des vari\'et\'es toriques
dans \cite{bourqui:eclate}, \cite{bourqui:deploye} et
\cite{bourqui:torique}.
\par
Ce texte est organis\'e de la fa\c con suivante:
dans la partie 1, nous fixons les notations et rappelons la d\'efinition
de la fonction z\^eta des hauteurs, dans la partie 2 nous g\'en\'eralisons
\`a la situation pr\'esente la d\'efinition de la mesure de Tamagawa
associ\'ee \`a une m\'etrique ad\'elique et
dans le paragraphe 3 nous \'enon\c cons le r\'esultat
dont la d\'emonstration occupe l'ensemble de la derni\`ere partie.
\par
Dans un souci de compl\'etude, nous redonnons l'ensemble des
constructions, bien que la plupart des notions soient
strictement analogues \`a celles d\'efinies sur un corps
de nombres. Les premi\`eres sections de ce texte contiennent
donc des redites par rapport \`a \cite{peyre:fano}.
\par\penalty-500
\section{Points de hauteur born\'ee}

Dans cette partie, nous transportons au cas d'un corps de fonctions global
la notion de syst\`eme de hauteurs connue sur un corps de nombres.

\subsection{Notations}
Dans la suite nous utiliserons les notations suivantes

\begin{notas}
Pour tout corps $\ff$, on note $\fbarre$ une cl\^oture
alg\'ebrique de $\ff$ et $\fs$ la cl\^oture s\'eparable de $\ff$
dans $\fbarre$.
\par
Si $\ff$ est un corps de fonctions global, on note $\Fqf$
le corps des constantes de $\ff$ et $\Mf$ l'ensemble
des places de $\ff$. Pour tout $\pp$ de $\Mf$, on note $\fp$
le compl\'et\'e de $\ff$ en $\pp$ et $\FF_\pp$ son corps r\'esiduel.
Pour toute place $\pp$ la norme $\norme$ sur $\fp$ est normalis\'ee
par la relation
\[\forall x\in\fp,\quad\normede x=(\cardinal\Fp)^{-\vp(x)}\]
o\`u $\vp:\fp\to\ZZ$ est la valuation surjective en $\pp$
et $\cardinal A$ d\'esigne le cardinal de $A$.
On dispose donc de la formule du produit
\[\forall x\in\ff,\quad\prod_{\pp\in\Mf}\normede x=1.\]
Pour toute place $\pp$, la mesure de Haar $\dxp$ sur $\fp$
est normalis\'ee par la relation
\[\int_{\Op}\dxp=1.\]
On note $\cc$ la courbe projective et lisse sur $\Fqf$
de corps de fonctions $\ff$, $\gf$ le genre de $\cc$
et $\hf$ le nombre de classes de diviseurs de degr\'e 0 sur $\cc$.
On identifiera $\Mf$ aux points ferm\'es de $\cc$.
\par
Nous omettrons $\ff$ dans ces notations lorsque le corps global 
sera indiqu\'e par le contexte.
\par
Si $\mV$ est un sch\'ema sur le spectre d'un anneau $A$
et si $B$ est une $A$ alg\`ebre commutative, alors
$\mV(B)$ d\'esigne l'ensemble
\[\Hom_{\Spec A}(\Spec B,\mV)\]
et $\mV_B$ le produit $\mV\times_{\Spec A}{\Spec B}$. Si $\VV$
est lisse sur un corps $\ff$, son groupe de Picard est
not\'e $\Pic\VV$, son groupe de N\'eron-Severi $\NS(\VV)$
et son faisceau canonique $\omegaV$.
Le c\^one de $\NS(\VV)\otimes_\ZZ\RR$ engendr\'e
par les classes de diviseurs effectifs est not\'e
$\CeffV$. On note \'egalement $\Vbarre$ la vari\'et\'e
$\VV_\fbarre$ et $\Vs=V_\fs$.
\par
Si $n$ n'est pas divisible par la caract\'eristique de $\ff$,
le faisceau \'etale sur $\VV$ des racines $n$-i\`emes de
l'unit\'e est not\'e $\mut n1$. Le faisceau constant $\mut n0$
est \'egalement not\'e $\mut n{{0}}$ et pour tout entier $j$
strictement positif, $\mut nj$ d\'esigne $\mut n{j-1}\otimes\mut n1$
et $\mut n{-j}$ le faisceau $\Hom(\mut nj,\mut n0)$.
Pour tout nombre premier $\Ll$ distinct de la caract\'eristique
de $\ff$
\[H\etale^i(\Vbarre,\ZZ_\Ll(j))=\inverselimit_n
H\etale^ i(\Vbarre,\mut{\Ll^n}j)\]
et
\[H\etale^i(\Vbarre,\QQ_\Ll(j))=
H\etale^ i(\Vbarre,\ZZ_\Ll(j))\otimes\QQ_\Ll\]
On note $\Br\VV$ le groupe de Brauer cohomologique de $\VV$ d\'efini par
\[\Br\VV=H\etale^2(V,\mathbf G_m).\]
\par
Si $\VV$ est une vari\'et\'e sur un corps global $\ff$, $\VV(\adeles)$
d\'esigne l'espace ad\'elique associ\'e tel qu'il est d\'efini par Weil dans
\cite[\S1]{weil:adeles}.
\end{notas}

\subsection{M\'etriques ad\'eliques et hauteurs}
La notion de m\'etrique ad\'elique est une g\'en\'eralisation 
imm\'ediate \`a notre cadre de celle d\'ecrite dans \cite{peyre:torseurs}
pour un corps de nombres, elle-m\^eme inspir\'ee
de la notion classique de hauteur (cf. \cite{lang:diophantine},
\cite{silverman:height}).

\begin{defi}
Soit $\VV$ une vari\'et\'e projective, lisse et g\'eom\'etriquement int\`egre
sur un corps global $\ff$ de caract\'eristique finie $\car$.
Soit $\LL$ un faisceau inversible sur $\VV$ et $\pp$ une place
de $\ff$. Une {\em m\'etrique $\pp$-adique} sur $\LL$ est une application
associant \`a tout point $x$ de $\VV(\fp)$ une norme $\metrique$
sur la fibre $\LL(x)=\LL_x\otimes_{\OVx}\fp$ de sorte que pour toute section
$\sectiondeL$ de $\LL$ d\'efinie sur un ouvert de $\VV$
l'application
\[x\mapsto\metriquede{\sectiondeL(x)}\]
soit continue pour la topologie $\pp$-adique.
\par
Si $\pp$ est une place de $\ff$ et $\mV$ un mod\`ele
projectif et lisse de $\VV$ sur $\Op$, alors \`a tout mod\`ele $\mL$
de $\LL$ sur $\mV$ on peut associer une m\'etrique $\pp$-adique
$\metrique$ sur $\LL$ qui est construite de la mani\`ere suivante:
$\VV$ \'etant projective, tout point $x$ de $\VV(\fp)$
d\'efinit un point $\tilde x$ de $\mV$ et $\tilde x^ *(\mL)$
d\'efinit une $\Op$-structure sur $\LL(x)$ dont on choisit un g\'en\'erateur
$y_0$; $\metrique$ est alors donn\'ee par la formule
\[\forall y\in\LL,\quad\metriquede y=\left\vert\frac{y}{y_0}\right\vert_\pp.\]
\par
Une famille de m\'etriques $\metriques$ sur $\LL$ est appel\'ee
{\em m\'etrique ad\'elique} s'il existe un sous-ensemble fini $\SSS$
de $\Mf$, un mod\`ele projectif et lisse $\mV$ de$\VV$
sur le compl\'ementaire de $\SSS$ dans $\cc$ et un mod\`ele $\mL$
de $\LL$ sur cet espace de sorte que pour toute place $\pp$
de $\Mf-\SSS$ la m\'etrique $\metrique$ soit d\'efinie par 
$\mL\otimes_{\structural_\mV}\Op$.
\par
Par abus de langage, nous appellerons {\em hauteur d'Arakelov} sur $\VV$
la donn\'ee d'une paire $\hauteur=(\LL,\metriques)$ o\`u $\LL$
est un faisceau inversible sur $\VV$ et $\metriques$
une m\'etrique ad\'elique sur ce fibr\'e. Pour tout point rationnel
$x$ de $\VV$ la {\em hauteur de $x$ relativement \`a $\hauteur$}
est alors d\'efinie par
\[\hauteur(x)=\prod_{\pp\in\Mf}\metriquede{\sectiondeL(x)}^{-1}\]
o\`u $\sectiondeL$ est une section de $\LL$ sur un voisinage de $x$,
non nulle en $x$.
\par
On note $\hauteurs(\VV)$ l'ensemble des classes d'isomorphisme 
de hauteurs d'Arakelov
modulo la relation d'\'equivalence $\sim$ engendr\'ee 
par les relations de la forme
\[(\LL,\metriques)\sim(\LL,(\lambda_\pp\metrique)_{\pp\in\Mf})\]
pour toute famille $(\lambda_\pp)_{\pp\in\Mf}$ de nombres r\'eels
strictement positifs v\'erifiant
$\prod_{\pp\in\Mf}\lambda_\pp=1$. Si $x$ est un point rationnel et
$\hauteur$ une hauteur d'Arakelov, $\hauteur(x)$ ne d\'epend que de
la classe de $\hauteur$ dans $\hauteurs(\VV)$.
\par
On dispose d'un morphisme d'oubli
\[\oubli:\hauteurs(\VV)\to \NS(\VV).\]
On appelle syst\`eme de hauteurs une section
$\HH$ de l'application d'oubli. Un syst\`eme de hauteurs $\HH$
sur $\VV$ induit un accouplement
\[\HH:\NS(\VV)\otimes\CC\times \VV(\ff)\to\CC\]
qui est l'exponentielle d'une fonction lin\'eaire en la premi\`ere variable
et telle que
\[\forall \LL\in\NS(\VV),\quad\forall x\in\VV(\ff),\quad
\HH(\LL,x)=\HH(\LL)(x).\]
\end{defi}

Nous renvoyons \`a \cite{peyre:cercle} pour des exemples
de hauteurs et de syst\`emes de hauteurs. En particulier pour tout morphisme
$\phi:\VV\to W$ de vari\'et\'es on a un diagramme commutatif
\[\begin{CD}
\hauteurs(W)@>\phi^*>>\hauteurs(\VV)\\
@VVV@VVV\\
\NS(W)@>\phi^*>>\NS(\VV)
\end{CD}\]
et si $\ee/\ff$ est une extension de corps on dispose d'un morphisme
de normes
\[\normeEaF:\hauteurs(\VV_\ee)\to\hauteurs(\VV).\]
Si $\HH_\ee$ est un syst\`eme de hauteurs sur $\VV_\ee$
la hauteur induite $\HH_\ff$ est d\'efinie par le diagramme
commutatif
\[\begin{CD}
\NS(\VV)@>>>\NS(\VV_\ee)\\
@VVN\HH_\ff V@VV\HH_\ee V\\
\hauteurs(\VV)@<\normeEaF<<\hauteurs(\VV_\ee)
\end{CD}\]
o\`u $N=[\ee:\ff]$.

\subsection{Sous-vari\'et\'es accumulatrices}

Comme dans le cas des corps de nombres deux notions 
de sous-vari\'et\'es accumulatrices
apparaissent naturellement.

\begin{defi}
Soit $\hauteur$ une hauteur d'Arakelov sur une vari\'et\'e
projective lisse et g\'eom\'etriquement int\`egre $\VV$
au-dessus d'un corps global $\ff$ de caract\'eristique $p$.
Pour tout sous-espace localement ferm\'e $W$ de $\VV$
et tout nombre r\'eel strictement positif $\unreel$, on pose
\[\nWH(\unreel)=\cardinal\{\,x\in W(\ff)\mid\hauteur(x)\varleq\unreel\,\}.\]
\end{defi}

\begin{rem}
si $\hauteur=(\LL,\metriques)$ avec $[\LL]$ appartenant \`a l'int\'erieur
de $\CeffV$, alors il existe un ouvert non vide $\UU$ de $\VV$
tel que $\nUH(\unreel)$ soit fini pour tout nombre r\'eel $\unreel$.
\end{rem}

\begin{defi}
On note pour tout sous-espace localement ferm\'e
$W$ de $\VV$
\[\aLde W=\mathop{\overline\lim}_{\unreel\to+\infty}\log
(\nWH(\unreel))/\log(\unreel)\varleq +\infty.\]
Un ferm\'e strict $\ferme$ de $\VV$ est dit $\LL$-accumulateur
au sens strict
si et seulement si pour tout ouvert non vide $W$ de $\ferme$,
il existe un ouvert non vide $\UU$ de $\VV$ avec
\[\aLde W>\aLde \UU\]
\par
Un ferm\'e strict $\ferme$ de $\VV$ est dit mod\'er\'ement accumulateur
pour $\hauteur$ si et seulement si, pour tout
ouvert non vide $W$ de $\ferme$, il existe un ouvert non vide $\UU$
de $\VV$ tel que
\[\mathop{\overline\lim}_{\unreel\to+\infty}\nWH(\unreel))/\nUH(\unreel)>0.\]
\end{defi}

\subsection{Fonction z\^eta des hauteurs}

Contrairement au cas des corps de nombres, les m\'etriques \'etant \`a
valeur dans $q^\ZZ$ sauf pour un nombre fini d'entre elles,
$\nUH(\unreel)$ n'est plus asymptotiquement \'equivalent \`a une fonction de
la forme $C\unreel^a(\log\unreel)^{b-1}$. Toutefois il est bien connu que
l'objet naturel dans le cadre fonctionnel est la fonction z\^eta associ\'ee,
et nous verrons plus loin qu'il est effectivement possible 
d'\'etendre \`a ce cadre les propri\'et\'es connues sur les corps de nombres.
De mani\`ere plus pr\'ecise, les travaux de Batyrev, Franke,
Manin et Tschinkel, (\cite[\S2]{fmt:fano}, \cite{batyrevtschinkel:toric}
et \cite{batyrevtschinkel:generaltoric}) ont soulign\'es
l'int\'er\^et de consid\'erer les fonctions z\^etas associ\'ees a un syst\`eme
de hauteurs et d\'efinies sur un ouvert du produit $\NS(\VV)\otimes\CC$.
Nous en rappelons maintenant la d\'efinition.

\begin{defis}
Soit $\HH$ un syst\`eme de hauteurs sur une vari\'et\'e
projective lisse et g\'eom\'etriquement int\`egre $\VV$
au-dessus d'un corps global $\ff$ de caract\'eristique $p$.
Soit $\UU$ un ouvert de $\VV$. La fonction z\^eta associ\'ee est d\'efinie
par la s\'erie
\begin{equation}
\label{equ:zeta:def}
\zetaHU(s)=\sum_{x\in\UU(\ff)}\HH(s,x)^{-1}
\end{equation}
o\`u $s$ d\'esigne un \'el\'ement de $\NS(\VV)\otimes\CC$
pour lequel la s\'erie converge absolument.
\end{defis}

Nous r\'eunissons maintenant quelques assertions sur le domaine
de convergence des fonctions z\^eta. Ces assertions sont bien connues dans
le cas de corps de nombres.

\begin{lemme}
L'ensemble sur lequel $\zetaHU$
converge absolument est un ensemble convexe de $\NS(\VV)\otimes\CC$
et si $s_0$ appartient \`a cet ensemble alors celui-ci contient
\[s_0+i\NS(\VV)\otimes\CC.\]
\end{lemme}

\begin{proof}
Ces assertions sont des propri\'et\'es g\'en\'erales des s\'eries z\^eta.
La convexit\'e r\'esulte de l'in\'egalit\'e de H\"older: en effet si
$\zetaHU$ converge absolument pour
\[s_0,s_1\in\NS(\VV)\otimes\CC\]
et si $\lambda_0,\lambda_1\in\RR_{>0}$ avec $\lambda_0+\lambda_1=1$,
alors pour tout sous-ensemble fini $I$ de $\UU(\ff)$, on a
\[\begin{split}
\sum_{x\in I}\vert\HH(\lambda_0s_0+\lambda_1s_1)\vert^{-1}&=
\sum_{x\in I}\vert\HH(s_0,x)\vert^{-\lambda_0}\vert
\HH(s_1,x)\vert^{-\lambda_1}\\
&\varleq(\sum_{i\in I}\vert\HH(s_0,x)\vert^{-1})^{\lambda_0}
\varleq(\sum_{i\in I}\vert\HH(s_1,x)\vert^{-1})^{\lambda_1}.
\end{split}\]
La seconde assertion r\'esulte du fait que pour tout $x$ de $\VV(\ff)$
on a
\[\forall s_0\in\NS(\VV)\otimes\CC,\quad
\forall s_1\in i\NS(\VV)\otimes\RR,\quad
\vert\HH(s_0+s_1,x)\vert=\vert\HH(s_0,x)\vert.\]
\end{proof}

\begin{lemme}
Si $\Pic(\Vbarre)$ est de rang fini et si $\CeffV$ est un c\^one
polyh\'edrique rationnel, c'est-\`a-dire s'il existe
$n_1,\dots,n_m$ dans $\NS(\VV)$ tels que
\[\CeffV=\sum_{i=1}^m\RR_{\vargeq 0}n_i\]
alors il existe un ouvert dense $\UU$ de $\VV$ tel que pour tout
$s$ de $\NS(\VV)\otimes\CC$ en lequel $\zetaHU$ converge absolument,
$s+\CeffV$ est contenu dans le domaine de convergence absolu
de $\zetaHU$.
\end{lemme}

\begin{proof}
Fixons des diviseurs effectifs $D_1,\dots,D_m$ tels que
\[\CeffV=\sum_{i=1}^m\RR_{\vargeq 0}[D_i]\]
et posons
\[\UU=\VV-\cup_{i=1}^m\Supp D_i.\]
Soit $\sectiondeL_i$ une section de $[D_i]$ correspondant \`a
$D_i$ pour $1\varleq i\varleq m$. Soient $s$ un \'el\'ement de
$\NS(\VV)\otimes\CC$ en lequel $\zetaHU$ converge absolument
et $s'$ un \'el\'ement de $\CeffV\cap\NS(\VV)\otimes\QQ$.
Il existe un entier $\lambda$ strictement positif tel que $\lambda s'$
puisse se mettre sous la forme
\[\lambda s'=\sum_{i=1}^mn_i[D_i]\]
avec $n_i\in\ZZ_{\vargeq 0}$ pour $1{\vargeq}i{\vargeq}m$.
Soit $(\LL,\metriques)$ un repr\'esentant de
$\HH(\sum_{i=1}^mn_i[D_i])$. Comme $\Pic(\overline\VV)$
est de rang fini, on peut supposer, quite \`a augmenter $\lambda$,
que la classe de $\LL$ co\"\i ncide avec celle
de $\sum_{i=1}^mn_i[D_i]$ dans $\Pic\VV$.
Le produit tensoriel $\bigotimes_{i=1}^{m}\sectiondeL_i^{n_i}$
d\'efinit alors une section $\sectiondeL$ de $\LL$.
Pour tout place $\pp$ de $\ff$, la fonction
\[\begin{array}{rcl}
\VV(\fp)&\to&\RR_{\vargeq 0}\\
x&\mapsto&\metriquede{\sectiondeL(x)}
\end{array}\]
est continue et admet un maximum $B_\pp$.
En outre si $\SSS$ est une partie finie de $\Mf$
et $\mV$ (resp. $\mL$) un mod\`ele de $\VV$ (resp. $\LL$)
sur $\cc-\SSS$ de sorte que les mestriques $\metrique$
soient d\'efinies par $\mL$ en-dehors de $\SSS$;
alors, quitte \`a augmenter $\SSS$, on peut
supposer que $\sectiondeL$ provient d'une section $\tilde\sectiondeL$
de $\mL$. Pour toute place $\pp$ dans $\Mf-\SSS$, tout point
$x$ de $\VV(\fp)$ correspondant \`a un point $\tilde x$
de $\mV(\Op)$, soit $y_0$ un g\'en\'erateur de $\tilde x^*(\mL)$.
On a alors la relation
\[\metriquede{\sectiondeL(x)}=\left\vert\frac{\sectiondeL(x)}{y_0}
\right\vert_\pp\]
mais, comme $\sectiondeL$ provient de $\tilde\sectiondeL$,
on a que $\sectiondeL(x)$ appartient \`a $u_0\Op$. Donc $B_\pp\varleq 1$.
\par
Si $\UU(\ff)=\emptyset$, alors le r\'esultat est \'evident.
Dans la cas contraire, soit $x\in\UU(\ff)$ pour
presque toute place $\pp$ de $\ff$ on a $\metriquede{\sectiondeL(x)}=1$,
par cons\'equent $B_\pp=1$ pour presque toute place $\pp$ de $\ff$. Mais alors
\[\forall x\in\UU(\ff),\quad
\HH(\lambda s',x)=
\prod_{\pp\in\Mf}\metriquede{\sectiondeL(x)}^{-1}\vargeq
\prod_{\pp\in\Mf}B_\pp^{-1}>0.\]
Autrement dit la fonction $x\mapsto\HH(\lambda s',x)$
est uniform\'ement minor\'ee sur $\UU(\ff)$. On obtient
\[\begin{split}
\sum_{x\in\UU(\ff)}\vert\HH(s+\lambda s')\vert^{-1}&
\varleq\sum_{x\in\UU(\ff)}\vert\HH(s,x)\vert^{-1}
\HH(\lambda s',x)^{-1}\\
&\vargeq\left(\prod_{\pp\in\Mf}B_\pp\right)\zetaHU(s).
\end{split}\]
La s\'erie d\'efinissant $\zetaHU$ converge donc absolument
en tout point de
\[s+\CeffV\cap\NS(\VV)\otimes\QQ\]
mais l'enveloppe convexe de ce c\^one est $s+\CeffV$
et le r\'esultat d\'ecoule du lemme pr\'ec\'edent.
\end{proof}

\begin{lemme}
Si $\LL$ est un faisceau inversible tr\`es ample sur $\VV$
et $N$ la dimension de $\Gamma(\VV,\LL)$, alors pour toute hauteur
$\hauteur$ de la forme $(\LL,\metriques)$ et tout $\epsilon>0$
la s\'erie $\sum_{x\in\VV(\ff)}\hauteur(x)^{-N-\epsilon}$
converge.
\end{lemme}

\begin{proof}
Soit $\sectiondeL_1,\dots,\sectiondeL_N$ une base de $\Gamma(\VV,\LL)$.
Le syst\`eme de m\'etriques d\'efinie par
\[\forall \pp\in\Mf,\quad\forall x\in\VV(\fp),\quad\forall y\in\LL(x),\quad
\metriquede y'=\sup_{
\atop{1\varleq i\varleq n}{\sectiondeL_i(x)\neq 0}}
\left\vert\frac{y}{\sectiondeL_i(x)}\right\vert_\pp\]
est une m\'etrique ad\'elique sur $\LL$. On note $\hauteur'$
la hauteur
$(\LL,(\metrique')_{\pp\in\Mf})$. Pour presque toute place $\pp$ de $\ff$,
$\metrique'$ co\"\i ncide avec $\metrique$. Il existe donc une constante $C>0$
telle que
\[\forall x\in\VV(\ff),\quad\hauteur(x)\vargeq C\hauteur'(x).\]
Il suffit donc de montrer le r\'esultat pour $\hauteur'$.
Mais
\[\sum_{x\in\VV(\ff)}\hauteur'(x)^{-N-\epsilon}\varleq
\sum_{x\in\PP^{N-1}(\ff)}\hauteur_{\PP^{N-1}}(x)^{-N-\epsilon}\]
o\`u $\hauteur_{\PP^{N-1}}$ est la hauteur usuelle sur $\PP^{N-1}$.
Or il r\'esulta de \cite[theorem, p. 19]{serre:mordellweil}
qu'il existe une constante $C$ telle que
\[\cardinal\{x\in\PP^{N-1}(\ff)\mid\hauteur_{\PP^{N-1}}(x)
=\qq^d\}<C(\qq^d)^N\]
Par cons\'equent
\[\sum_{x\in\PP^{N-1}(\ff)}\hauteur_{\PP^{N-1}}(x)^{-N-\epsilon}
\varleq C\sum_d(\qq^d)^{-N-\epsilon}\qq^{dN}\qed\]
\noqed
\end{proof}

Pour terminer, nous allons \'enoncer une condition
qui implique trivialement la p\'eriodicit\'e
la fonction z\^eta des hauteurs

\begin{defi}
Dans la suite, nous dirons que le syst\`eme de hauteurs $\HH$
v\'erifie la condition (P) si et seulement si pour tout
\'el\'ement $\LL$ de $\NS(\VV)$, il existe un repr\'esentant
$(\LL,\metriques)$ de $\HH(\LL)$ dont les m\'etriques sont toutes \`a valeurs dans $\qq^{\ZZ}$.
\end{defi}

\begin{rem}
Il r\'esulte des d\'efinitions que
les fonctions $\zetaHU$ sont alors p\'eriodiques de groupe de p\'eriodes
contenant
$2i\pi/\log(\qq)\NS(\VV)$.
\end{rem}
\section{Mesures de Tamagawa}
\subsection{Quelques hypoth\`eses sur les vari\'et\'es}
\label{subsection:mesures:hypos}

Dans ce paragraphe, nous allons pr\'eciser quelques hypoth\`eses sur
les vari\'et\'es qui nous permettront de d\'efinir la mesure de Tamagawa
associ\'ee \`a une m\'etrique ad\'elique sur le faisceau anticanonique.
Les hypoth\`eses g\'eom\'etriques sont automatiquement v\'erifi\'ees
par une vari\'et\'e de Fano sur un corps de caract\'eristique $0$.
En l'absence d'un th\'eor\`eme d'annulation de Kodaira cela n'est
toutefois plus le cas en caract\'eristique finie.

\begin{hyposgeo}
Dans la suite $\VV$ d\'esigne une vari\'et\'e projective
lisse et g\'eom\'etriquement int\`egre sur un corps global $\ff$
de caract\'eristique $\car$ finie
v\'erifiant les assertions suivantes:
\begin{enumerer}
\item la classe du faisceau anticanonique $\antican$
dans $\NS(\VV)\otimes\RR$ appartient \`a l'int\'erieur du c\^one
$\CeffV$,
\item les groupes de cohomologie $H^i(V,\OV)$
sont nuls pour $i=1$ ou $2$,
\item \label{item:hypo:picard}
le groupe $\Pic(\Vs)$ est un $\ZZ$-module libre de rang fini
et co\"\i ncide avec $\Pic\Vbarre$,
\item le c\^one $\CeffVbarre$ est poly\'edrique rationnel, c'est-\`a-dire
qu'il existe $m_1,\dots,m_r$ dans $\Pic(\Vbarre)$
tels que
\[\CeffVbarre=\sum_{i=1}^rm_i,\]
\item \label{item:hypo:brauer}
si $\Ll$ est un nombre premier distinct de $\car$, la partie $\Ll$-primaire
de $\Br\Vbarre$ est finie.
\end{enumerer}
\end{hyposgeo}
\begin{rem}
Si $\VV$ est $\fs$-rationnelle alors la premi\`ere assertion de 
(\ref{item:hypo:picard}) et (\ref{item:hypo:brauer})
sont v\'erifi\'ees. En effet, par \cite[appendice 2.A]{cts:descente2},
si $\VV$ est $\fs$-rationnelle alors $\Pic(\Vs)$ est stablement
isomorphe \`a $\ZZ$ et donc libre de rang fini sur $\ZZ$
et par \cite[corollaire 7.5]{grothendieck:brauerIII},
la partie $\Ll$-primaire du groupe de Brauer esy un invariant birationnel
des vari\'et\'es propres et lisses, ce qui implique la trivialit\'e
de ce groupe si $\VV$ est $\fs$-rationnelle.
\end{rem}

Dans la suite nous supposerons \'egalement que la vari\'et\'e $\VV$ v\'erifie
la condition suivante:

\begin{hypoari}
L'ensemble $\VV(\ff)$ des points rationnels de $\VV$ est
dense pour la topologie de Zariski.
\end{hypoari}

\subsection{Ensembles de mauvaises places}

Comme dans le cas des corps de nombres, les facteurs
de convergence qui apparaissent naturellement sont les facteurs
locaux de la fonction $L$ associ\'ee au groupe de Picard
g\'eom\'etrique de la vari\'et\'e. Le choix de ces facteurs
est justifi\'e par les arguments de descente qui les font appara\^\i tre
dans le passage aux torseurs universels (cf. \cite{salberger:tamagawa}
et \cite{peyre:torseurs}). Mais pour montrer la convergence des mesures de
Tamagawa, il faut comparer ces facteurs locaux \`a ceux donn\'es par
le deuxi\`eme groupe de cohomologie $\Ll$-adique,
ce qui est rendu possible par le lemme suivant:

\begin{lemme}
\label{lemme:mauvaises:places}
Sous les hypoth\`eses du paragraphe \ref{subsection:mesures:hypos},
il existe un ensemble fini de places $\SSS$ et un mod\`ele projectif
et lisse $\mV$ de $\VV$ sur $\cc-\SSS$ dont les fibres sont 
g\'eom\'etriquement int\`egres et tel que pour toute place $\pp$
en-dehors de $\SSS$,
\begin{itemize}
\item Il existe un isomorphisme de $\Pic(\Vbarre)$
sur $\Pic(\mV_{\Fpbarre})$ compatible avec les actions des groupes de Galois,
\item Pour tout nombre premier $\Ll$ distinct de $\car$,
la partie $\Ll$-primaire du groupe
$\Br(\mV_\Fpbarre)$ soit finie.
\end{itemize}
\end{lemme}

\begin{proof}
La vari\'et\'e $\VV$ \'etant projective, on la plonge dans un 
espace projectif $\PP^N_\ff$. Soit $\mV$
l'adh\'erence de $\VV$ dans $\PP^N_{\cc}$. Comme $\VV$
est lisse et g\'eom\'etriquement int\`egre, il existe
par \cite[IV 6.8.7 et 9.7.7]{egan} un ensemble fini $\SSS$
de places tel que $\mV\times_{\cc}\cc-\SSS$ soit lisse \`a
fibres g\'eom\'etriquement int\`egres sur $\cc-\SSS$.
\par
Par hypoth\`ese, $\Pic(\Vbarre)$ est un $\ZZ$-module libre
de type fini qui co\"\i ncide avec $\Pic(\Vs)$. Il existe
donc une extension galoisienne finie $\ee$ de $\ff$
telle que $\VV(\ee)\neq\emptyset$ et
\[\Pic(\VV_\ee)\iso\Pic(\Vbarre).\]
On ajoute \`a $\SSS$ les places qui se ramifient
dans l'extension $\ee/\ff$ et on note $\SSS_\ee$ l'ensemble
des places de $\ee$ au-dessus de $\SSS$. Soit $\pP$ une place
de $\ee$ en-dehors de $\SSS_\ee$, soit
$\widehat{\OP\hs}$ le compl\'et\'e d'un hens\'elis\'e strict de
$\OP$ et $\Fr(\widehat{\OP\hs})$ son corps des fractions.
Le sch\'ema $\mV_{\widehat{\OP\hs}}$ \'etant lisse,
l'application canonique
\[\Pic(\mV_{\widehat{\OP\hs}})\to\Pic(\VV_{\Fr(\widehat{\OP\hs})})\]
est un isomorphisme. Comme, par hypoth\`ese, $\Pic(\VV_\ee)$
est isomorphe \`a $\Pic(\Vbarre)$ et donc \`a $\NS(\Vbarre)$
et que le groupe de N\'eron-Severi ne change pas par extension de
corps alg\'ebriquement clos, l'application
\[\Pic(\VV_\ee)\to\Pic(\VV_{\Fr(\widehat{\OP\hs})})\]
est \'egalement un isomorphisme. Consid\'erons alors
la compos\'ee des morphisme naturels
\begin{equation}
\label{equ:mesures:picard}
\Pic(\VV_\ee)\iso\Pic(\VV_{\Fr(\widehat{\OP\hs})})\iso
\Pic(\mV_{\widehat{\OP\hs}})\to\Pic(\mV_\FPbarre)
\end{equation}
Par hypoth\`ese, $H^i(\VV,\OV)=0$ pour $i=1$ ou $2$.
Par le th\'eor\`eme de semi-continuit\'e (cf. 
\cite[theorem III.12.8]{hartshorne:geometry}), on peut,
quitte \`a augmenter $\SSS$, supposer que
$H^i(\mV_{\Fp},\anneau_{\mV_{\Fp}})$ est nul pour $i=1$
ou $2$ et $\pp\in\Mf-\SSS$. Il en r\'esulte que
$H^i(\mV_\Fpbarre,\anneau_{\mV_\Fpbarre})$ est trivial
pour $i=1$ ou $2$.
Par \cite[corollaire 1 de la proposition 3]{grothendieck:formel}
la fl\`eche de droite dans \eqref{equ:mesures:picard}
est \'egalement un isomorphisme. On obtient ainsi
l'isomorphisme recherch\'e.
\par
En ce qui concerne la seconde assertion, le corang $\Ll$-adique de 
$\Br(\mV_\Fpbarre)$ co\"\i ncide, d'apr\`es 
\cite[corollaire 3.4]{grothendieck:brauerII}
avec la diff\'erence entre $B_2$, le deuxi\`eme nombre de Betti
$\Ll$-adique de $\mV_\Fpbarre$ et $\rho$ le rang du groupe de Picard
de $\mV_\Fpbarre$. Par ce qui pr\'ec\`ede $\rho$ est
\'egalement le rang de $\Pic\Vbarre$ et par \cite[page 19-02]{serre:facteurs}
$B_2$ co\"\i ncide avec le deuxi\`eme nombre de Betti de $\Vbarre$.
Par hypoth\`ese la partie $\Ll$-primaire de $\Br(\Vbarre)$ est finie
et $B_2=\rho$.
\end{proof}

\begin{conv}
Dans la suite de ce texte, les paires $(\SSS,\mV)$
consid\'er\'ees avec $\SSS$ un ensemble fini de places de $\ff$
et $\mV$ un mod\`ele de $\VV$ sur $\cc-\SSS$ v\'erifient les conditions
du lemme.
\end{conv}

\subsection{Mesures locales}
Dans ce paragraphe, nous associons \`a toute m\'etrique
ad\'elique $\metriques$ sur $\antican$ des mesures $\mesp$
sur $\VV(\fp)$.

\begin{defi}
Soit $\pp$ une place de $\ff$ et $\metrique$
une m\'etrique $\pp$-adique sur $\antican$. Si $x$
appartient \`a $\VV(\fp)$, on choisit des coordonn\'ees
locales $x_1,\dots,x_n$ en $x$ d\'efinissant un morphisme de
$\ff$ vari\'et\'es $f$ d'un ouvert de Zariski $\UU$ de $\VV$
dans $\affine^n_\ff$ et induisant un isomorphisme analytique
pour la topologie $\pp$-adique
d'un ouvert $W$ de $V(\fp)$ sur son image.
Le morphisme de faisceau coh\'erent
\[\omega(f):f^*(\canoniquede{\affine^n_\ff/\ff})
\to\canoniquede{U/\ff}\]
induit un isomorphisme de faisceau
pour la topologie $\pp$-adique
\[{}^t\omega(f)^{-1}:f^*\canoniquede{f(W)}^{-1}\to
\canoniquede W^{-1}\]
qui permet d'identifier 
$\frac{\partial}{\partial x_1}\wedge\dots\wedge
\frac{\partial}{\partial x_n}$ avec une section de
$\canoniquede W^{-1}$. Sur $W$, la mesure est alors d\'efinie
par la relation
\[\mesp=\left\Vert\frac{\partial}{\partial x_1}\wedge\dots\wedge
\frac{\partial}{\partial x_n}\right\Vert_\pp
\Haar{x_{1,\pp}}\dots\Haar{x_{n,\pp}}\]
o\`u $\Haar{x_{i,\pp}}$ d\'esigne la mesures de Haar 
normalis\'ee sur $\fp$.
\par
Si $x_1,\dots,x_n$ et $x'_1,\dots,x'_n$
sont deux syst\`emes de coordonn\'ees d\'efinies sur un m\^eme ouvert
$W$ et correspondant \`a des fonctions
\[f,f':\UU\to\affine^n_\ff\]
Soit $\Phi$ l'isomorphisme analytique
\[f'\circ f^{-1}:f(W)\to f'(W)\]
On a alors la relation
\[{}^t\omega(\Phi)^{-1}\left(\frac{\partial}{\partial x'_1}\wedge\dots\wedge
\frac{\partial}{\partial x'_n}\right)=\Jac_x(\Phi)^{-1}
\frac{\partial}{\partial x_1}\wedge\dots\wedge
\frac{\partial}{\partial x_n}\]
et par \cite[\S2.2.1]{weil:adeles}, on a
\[\Haar{x'_{1,\pp}}\dots\Haar{x'_{n,\pp}}=
\normede{\Jac_x(\Phi)}\Haar{x_{1,\pp}}\dots\Haar{x_{n,\pp}}\]
On en d\'eduit les \'egalit\'es
\[\begin{split}
&\left\Vert\frac{\partial}{\partial x'_1}\wedge\dots\wedge
\frac{\partial}{\partial x'_n}\right\Vert_\pp
\Haar{x'_{1,\pp}}\dots\Haar{x'_{n,\pp}}\\
&=\left\Vert\frac{\partial}{\partial x_1}\wedge\dots\wedge
\frac{\partial}{\partial x_n}\right\Vert_\pp
\normede{\Jac_x(\Phi)}^{-1}
\Haar{x'_{1,\pp}}\dots\Haar{x'_{n,\pp}}\\
&=\left\Vert\frac{\partial}{\partial x_1}\wedge\dots\wedge
\frac{\partial}{\partial x_n}\right\Vert_\pp
\Haar{x_{1,\pp}}\dots\Haar{x_{n,\pp}}.
\end{split}\]
Les mesures obtenues se recollent donc pour former une mesure bor\'elienne
$\mesp$ sur $\VV(\fp)$.
\end{defi}

\subsection{Lien avec les densit\'es locales}

Dans ce paragraphe, nous allons faire le lien pour $\pp\in\Mf-\SSS$
entre le volume $\mesp(V(\fp))$ et la densit\'e locale en $\pp$
d\'efinie par
\[d_\pp(\VV)=\frac{\cardinal\mV(\Fp)}{(\cardinal\Fp)^{\dim\VV}}.\]
La d\'emonstration suit celle de \cite[\S2.2.1]{weil:adeles},
\cite[\S2.2.2]{peyre:fano} et \cite[theorem 2.13]{salberger:tamagawa}.

\begin{lemme}
Pour presque toute place $\pp$ de $\Mf-\SSS$, on a:
\[\mesp(\VV(\fp))=d_\pp(\VV).\]
\end{lemme}

\begin{proof}
Soit $n$ la dimension de $\VV$.
Fixons un plongement $\Phi:\VV\to\PP^N_\ff$ et $\mV$
l'adh\'erence de $\VV$ dans $\PP^N_{\cc}$. Par le crit\`ere
valuatif de propret\'e, on a une bijection de $\VV(\fp)$
sur $\mV(\Op)$ qui induit pour tout $m$ des applications
surjectives
\[\pi_m:\VV(\fp)\to\mV(\Op/\pp^m).\]
On note $[x]_m=\pi_m^{-1}(\pi_m(x))$ pour tout $x$ de
$\VV(\fp)$. Soit $I$ le faisceau d'id\'eaux d\'efini par $\VV$
et $\calli I$ celui d\'efini par $\mV$. En dehors d'un nombre fini
de places, le sch\'ema $\mV_{\Op}$ est lisse et le faisceau
$(\Lambda^{n}\Omega^1_{\mV_{\Op}/\Op})\dual$ est un mod\`ele
de $\antican$. On peut donc supposer que la m\'etrique
$\metrique$ est d\'efinie par ce mod\`ele. Par
\cite[theorem 8.1, theorem 8.13]{hartshorne:geometry},
on a des suites exactes de faisceaux
\[\calli I/\calli I^2\stackrel\delta\to\Omega^1_{\PP^N_{\cc}/\cc}\otimes
\anneau_\mV\stackrel r\to\Omega^1_{\mV/\cc}\to 0\]
et
\[0\to J/J^2\stackrel\delta\to\Omega^1_{\PP^N_{\ff}/\ff}\otimes
\anneau_\VV\stackrel r\to\Omega^1_{\VV/\ff}\to 0.\]
Sur l'ouvert $\mV_i$ d\'efini par $X_i\neq 0$, on obtient des suites exactes
\begin{equation}
\label{equ:diff:un}
\calli I/\calli I^2_{\vert\mV_i}\stackrel\delta\to
\bigoplus_{j\neq i}\anneau_{\mV_i}\Haar X_j
\stackrel {r_i}\to\Omega^1_{\mV/\cc}\to 0
\end{equation}
et en notant $\VV_i={\mV_i}_\ff$ qu'on identifie
avec son image $U_i$ dans $\affine^{N+1}_\ff$,
\begin{equation}
\label{equ:diff:deux}
J/J^2_{\vert U_i}\stackrel\delta\to
\bigoplus_{j\neq i}\anneau_{U_i}\Haar X_j
\stackrel {r_i}\to\Omega^1_{U_i/\ff}\to 0
\end{equation}
Par cons\'equent pour presque toute place $\pp$ de $\Mf$, on a
\begin{equation}
\label{equa:norme:can}
\forall x\in U_i(\fp),\quad\forall y\in\antican(x),\quad
\metriquede y=\inf_{\quad\lower4pt\clap%
{$\scriptstyle\atop{0\leq i_1<\dots<i_n}{j\not\in\{i_1,\dots,i_n\}}$}\quad}
\normede{y(r_i(\Haar X_{i_1})\wedge\dots\wedge r_i(dx_{i_n}))}.
\end{equation}
Notons $f_i$ l'isomorphisme de $\VV_i$ sur $U_i$. On peut en outre
supposer que $\pp$ v\'erifie les conditions I et II de 
\cite[page 19]{weil:adeles} pour la famille $(f_i)_{0\leq i\leq N+1}$.
Soit $x\in\VV(\fp)$ et $(x_0:\dots:x_N)$ des coordonn\'ees homog\`enes
pour $\Phi(x)$. Apr\`es permutation des coordonn\'ees et 
multiplication par un scalaire, on peut supposer que $x_0=1$
et $x_1,\dots,x_n\in\Op$. Le point $x$ provient alors d'un \'el\'ement 
$\tilde x$ de $\calli U_0(\Op)$. Par \cite[theorem 2.2.3]{weil:adeles},
l'ensemble $[x]_1$ co\"\i ncide avec
\[\{(1:x'_1:\dots:x'_n)\in\VV(\fp)\mid x'_i\in x_i+\pp\}\]
On peut, \`a permutation  des coordonn\'ees pr\`es, supposer
par \eqref{equ:diff:un} que $\Omega^1_{\mV_{\Op},x}$
est isomorphe \`a $\bigoplus_{i=1}^n\anneau_{\mV_{\Op,x}}\Haar X_i$.
La famille $(X_1,\dots,X_n)$ constitue alors un syst\`eme de coordonn\'ees
locales au voisinage de $x$. Par \cite[page 22]{weil:adeles}
ce syst\`eme s'\'etend \`a $[x]_1$ et induit un diff\'eomorphisme
de $[x]_1$ sur $(x_1,\dots,x_n)+\pp^n$ et
l'isomorphisme entre $\Omega^1_{\mV_{\Op}/\Op,x}$
et $\bigoplus_{i=1}^n\anneau_{\mV_{\Op},x}\Haar X_i$
s'\'etend \'egalement \`a $[x]_1$. Par cons\'equent 
\eqref{equa:norme:can} se r\'e\'ecrit alors
\[\forall x'\in[x]_1,\quad\forall y'\in\antican(x'),\quad
\metriquede{y'}=\left\vert\frac{y'}{\frac{\partial}{\partial X_1}\wedge
\dots\wedge\frac\partial{\partial X_n}}\right\vert_\pp\]
et sur $[x]_1$, la mesure $\mesp$ co\"\i ncide avec
\[\Haar X_{1,\pp}\dots\Haar X_{n,\pp}.\]
On obtient donc
\[\mesp([x]_1)=\int_{(x_1,\dots,x_n)+\pp^n}\Haar X_{1,\pp}\dots
\Haar X_{n,\pp}=\cardinal\Fp^{-n}.\]
En somman sur $\mV(\fp)$, on obtient la formule du lemme.
\end{proof}

\subsection{Facteurs de convergence}
Comme dans \cite[\S2.2.3]{peyre:fano}, nous allons maintenant
appliquer les conjectures de Weil d\'emontr\'ees par Deligne
pour obtenir des facteurs de convergence.

\begin{defi}
Pour tout $\pp\in\Mf-\SSS$, on note $\Fr_\pp$ le morphisme de
Frobenius sur $\Fpbarre$ d\'efini par $x\mapsto x^{\cardinal\Fp}$.
La suite exacte
\[0\to I_\pp\to\Gal(\fpbarre/\fp)\to\Gal(\Fpbarre/\Fp)\to 0\]
o\`u $I_\pp$ d\'esigne le groupe d'inertie en $\pp$ et l'inclusion 
canonique du groupe $\Gal(\fpbarre/\fp)$ dans $\Gal(\fbarre/\ff)$
d\'efinissent une action de $\Gal(\Fpbarre/\Fp)$
sur $(\Pic(\Vbarre))^{I_\pp}$. On note \'egalement
$\Fr_\pp$ le morphisme de Frobenius induit sur $(\Pic(\Vbarre))^{I_\pp}$.
Le terme local de la fonction $L$ associ\'ee \`a $\Pic(\Vbarre)$
est alors d\'efini par
\[L_\pp(s,\Pic(\Vbarre))=
\frac{1}{\Det(1-\cardinal\Fp^{-s}\Fr_\pp\mid
(\Pic(\Vbarre))^{I_\pp}\otimes\QQ)}\]
\end{defi}

\begin{prop}
Pour toute place $\pp$ de $\Mf-\SSS$, le terme $L_\pp(1,\Pic(\Vbarre))$
est bien d\'efini et le produit
\[\prod_{\pp\in\Mf}\frac{d_\pp(\VV)}{L_\pp(1,\Pic(\Vbarre))}\]
est absolument convergent.
\end{prop}
\begin{proof}
Par la d\'emonstration du lemme \ref{lemme:mauvaises:places},
il existe une extension galoisienne finie $\ee$ de $\ff$
telle que $\Pic(\VV_\ee)\iso\Pic(\Vbarre)$. L'action
de $\Gal(\fbarre/\ff)$ sur $\Pic(\Vbarre)$ se factorise
donc par $\Gal(\ee/\ff)$ et pour tout
$\pp\in\Mf$,  on a que $(\Fr_\pp)^{[\ee:\ff]}$
agit trivialement sur
$\Pic(\Vbarre)^{I_\pp}$. Les valeurs propres de $\Fr_\pp$
sont donc des racines $[\ee:\ff]$-i\`emes de l'unit\'e
et
\[\Det(1-\cardinal\Fp^{-s}\Fr_\pp\mid(\Pic(\Vbarre)^{I_\pp}\otimes\QQ)\]
est non nul ce qui montre la premi\`ere assertion.
\par
Par la formule de Lefschetz (cf. \cite{serre:frobenius}) on a pour
toute place $\pp$ en-dehors de $\SSS$ et tout nombre premier $\Ll$
distinct de $p$,
\[\cardinal\mV(\Fp)=\sum_{i=0}^{2n}(-1)^i\Tr(\Fr_\pp\mid
H^i\etale(\mV_{\Fpbarre},\QQ_\Ll))\]
o\`u $n$ d\'esigne la dimension de $\VV$ et $\Fr_\pp$
le Frobenius g\'eom\'etrique. La vari\'et\'e $\mV_{\Fpbarre}$
\'etant lisse, projective et g\'eom\'etriquement int\`egre,
on a un isomorphisme canonique
\[H^{2n}\etale(\mV_{\Fpbarre},\QQ_\Ll(n))\iso\QQ_\Ll.\]
D'autre part les suites exactes de Kummer
\[0\to\mut{\Ll^r}1\to\mathbf G_m@>\times\Ll^r>>\mathbf G_m\to 0\]
o\`u $r$ est positif induisent des suites exactes
\[0\to H^1\etale(\mV_{\Fpbarre},\ZZ_\Ll(1))
\to T_\Ll(\Pic(\mV_{\Fpbarre}))\to 0\]
et
\[0\to\Pic(\mV_{\Fpbarre})\otimes\ZZ_\Ll\to
H^2\etale(\mV_{\Fpbarre},\ZZ_\Ll(1))\to
T_\Ll(\Br(\mV_{\Fpbarre}))\to0.\]
Mais il r\'esulte du lemme \ref{lemme:mauvaises:places}
que les modules de Tate
\[T_\Ll(\Pic(\mV_{\Fpbarre}))=\inverselimit_n\Pic(\mV_{\Fpbarre})_{(\Ll^n)}\]
et 
\[T_\Ll(\Br(\mV_{\Fpbarre}))=\inverselimit_n\Br(\mV_{\Fpbarre})_{(\Ll^n)}\]
sont tous les deux nuls. On obtient donc la trivialit\'e du groupe
$H^1\etale(\mV_{\Fpbarre},\QQ_\Ll(1))$ et un isomorphisme
\[H^2\etale(\mV_{\Fpbarre},\QQ_\Ll(1))\iso
\Pic(\mV_{\Fpbarre})\otimes\QQ_\Ll\]
qui, par dualit\'e de Poincar\'e, entra\^\i ne la
trivialit\'e de $H^{2n-1}\etale(\mV_{\Fpbarre},\QQ_\Ll)$
et induit un isomorphisme
\[H^{2n-2}\etale(\mV_{\Fpbarre},\QQ_\Ll(n-1))\iso
(\Pic(\mV_{\Fpbarre})\otimes\QQ_\Ll)\dual.\]
Or pour toute paire d'eniers $i$ et $j$, on a
\[\Tr(\Fr_\pp\mid H^i\etale(\mV_{\FPbarre},\QQ_\Ll(j))
=\frac{1}{\cardinal\Fp^j}\Tr(\Fr_\pp\mid H^i\etale(\mV_{\Fpbarre},\QQ_\Ll).\]
Par cons\'equent
\[d_\pp(\VV)=1+\frac{1}{\cardinal\Fp}\Tr(\Fr_\pp\mid
\Pic(\mV_{\Fpbarre})\otimes\QQ_\Ll)+
\sum_{i=0}^{2n-3}\frac{(-1)^i}{\cardinal\Fp^{\dim\VV}}
\Tr(\Fr_\pp\mid H^i\etale(\mV_{\Fpbarre},\QQ_\Ll)).\]
Mais par la conjecture de Weil d\'emontr\'ee par Deligne 
sur les valeurs propres des endomorphismes de Frobenius 
\cite[theorem 1.6]{deligne:weil}, on a des in\'egalit\'es
\[\vert\Tr(\Fr_\pp\mid H^i\etale(\mV_{\Fpbarre},\QQ_\Ll))\vert\leq
\cardinal\Fp^{i/2}\dim H^i\etale(\mV_{\Fpbarre},\QQ_\Ll).\]
Or le $i$-i\`eme nombre de Betti \'etale 
$\dim H^i\etale(\mV_{\Fpbarre},\QQ_\Ll)$ est constant sur les places
$\pp$ de bonne r\'eduction (cf. \cite[page 19-02]{serre:facteurs})
et donc
\[d_\pp(\VV)=1+\frac{1}{\cardinal\Fp}\Tr(\Fr_\pp\mid\Pic(\mV_{\Fpbarre})
\otimes\QQ)+O\left(\frac{1}{\cardinal\Fp^{3/2}}\right)\]
D'autre part les valeurs propres de $\Fr_\pp$
sur $\Pic(\mV_{\Fpbarre})$ qui est isomorphe
\`a $\Pic(\VV_\ee)$ \'etant des racines de l'unit\'e, on a \'egalement
\[\Det(1_\cardinal\Fp^{-1}\Fr_\pp\mid
\Pic(\mV_{\Fpbarre})\otimes\QQ)=
1+\frac{1}{\cardinal\Fp}\Tr(\Fr_\pp\mid\Pic(\mV_{\Fpbarre})\otimes\QQ)
+O\left(\frac{1}{\cardinal\Fp^{3/2}}\right)\]
et il en r\'esulte que
\[\frac{d_\pp(V)}{L\pp(1,\Pic(\Vbarre))}=1+O(\frac{1}{\cardinal\Fp^{3/2}}\]
et le produit de ces termes converge absolument.
\end{proof}

\begin{defi}
Pour toute place $\pp$ de $\ff$, on pose
\[\coeffconv=L_\pp(1,\Pic(\Vbarre)).\]
\end{defi}

\subsection{Mesure de Tamagawa}
Afin de normaliser le mesure nous aurons besoin 
de prendre le r\'esidu de la fonction $L$
associ\'ee \`a $\Pic\Vbarre$
et donc du lemme suivant:

\begin{lemme}
Le produit eul\'erien
\[\prod_{\pp\in\Mf}L_\pp(s,\Pic(\Vbarre))\]
converge absolument pour $\reel s>1$ et d\'efinit une fonction
$L(s,\Pic(\Vbarre))$ qui se prolonge en une fonction rationnelle de $q^{-s}$
sur $\CC$ avec un p\^ole d'ordre $t=\rg\Pic(\VV)$ en $s=1$.
\end{lemme}

\begin{proof}
La convergence pour $\reel s>1$ r\'esulte de la d\'efinition
et du fait que les valeurs propres du Frobenius agissant sur
$\Pic(\Vbarre)$ sont des racines de l'unit\'e.
\par
Soit $\Fr_{\qf}$ l'\'el\'ement du groupe $\Gal(\fbarre/\ff)$
envoyant $x$ sur $x^{\qf}$. La fonction $L$ est alors 
donn\'ee par le produit eul\'erien
\[L(s,\Pic(\Vbarre)=\prod_{x\in{\cc}_{(0)}}\frac{1}{\Det(1-\cardinal
\kappa(x)^{-s}\Fr_{\qf}^{[\kappa(x):\Fqf]}\mid\Pic(\Vbarre)\otimes\QQ)}\]
o\`u ${\cc}_{(0)}$ d\'esigne l'ensemble des points ferm\'es de $\cc$
et $\kappa(x)$ le corps r\'esiduel en $x$. Notons
$(\alpha_i)_{i\in I}$ les valeurs propres de $\Fr_{\qf}$ agissant
sur $\Pic(\Vbarre)\otimes\QQ$
et $(m_i)_{i\in I}$ leurs multiplicit\'es. On obtient
\[\begin{split}
L(s,\Pic(\Vbarre))&=\prod_{i\in I}\prod_{x\in{\cc}_{(0)}}
\frac{1}{(1-(q^{-s}\alpha_i)^{[\kappa(x):\Fqf]})^{m_i}}\\
&=\prod_{i\in I}Z(\cc,q^{-s}\alpha_i)
\end{split}\]
o\`u $Z(\cc,t)$ est la fonction z\^eta de $\cc$ d\'efinie
par
\[Z(\cc,t)=\exp\left(\sum_{n>0}\frac{\cardinal
\cc(\mathbf F_{\qf^n})t^n}{n}\right).\]
Par le th\'eor\`eme de Weil $Z(\cc,t)$ est une fonction rationnelle
de $t$ avec un p\^ole d'odre $1$ en $t=q^{-1}$. De mani\`ere
plus pr\'ecise, on a en fait
\[Z(\cc,t)=\frac{\Det(1-t\Fr_{\qf}\mid 
H^1\etale(\cc_{\overline{\mathbf F}_{\qf}}))}{(1-t)(1-qt)}\]
Ceci implique la deuxi\`eme assertion.
\end{proof}
\begin{defi}
On consid\`ere la mesure de Tamagawa
\[\mesH=\left(\lim_{s\to 1}(s-1)^tL(s,\Pic(\Vbarre))\right)
\frac{1}{\qf^{(\gf-1)\dim\VV}}\prod_{\pp\in\Mf}\coeffconv^{-1}\mesp\]
et le nombre de Tamagawa de $\VV$ relativement \`a la hauteur $\hauteur$
est d\'efini par
\[\tauV=\mesH(\overline{\VV(\ff)})\]
o\`u $\overline{\VV(\ff)}$ d\'esigne l'adh\'erece des
points rationnels de $\VV$ dans l'espace ad\'elique $\VV(\adeles)$.
\end{defi}
\section{Pr\'esentation du r\'esultat}
\subsection{Facteurs $\alphaV$ et $\betaV$}

Nous allons maintenant d\'efinir l'analogue de la fonction
caract\'eristique $\chiCeff$, qui remonte \`a
\cite{kocher:bereiche}, et qui fut introduite dans le contexte 
des conjectures de Manin par Batyrev et Tschinkel dans
\cite{batyrevtschinkel:toric}. Cette analogue est fourni par le facteur local
de la fonction $L$ de Draxl associ\'ee au c\^one effectif
(cf. \cite{draxl:tori}).

\begin{defi}
Soit $M$ un $\ZZ$-module libre et $C$ un c\^one rationnel 
poly\'edrique strictement convexe de $M\otimes\RR$; 
c'est-\`a-dire qu'il existe une famille finie
$(m_i)_{1\varleq i\varleq r}$ d'\'el\'ements de $M$
telle que $C=\sum_{i=1}^r\Rplus m_i$ avec
$C\cap -C=\{0\}$. On note $C\dual$
le c\^one dual d\'efini par
\[C\dual=\{\,y\in(M\otimes\RR)\dual\mid\forall x\in C,\,
\langle x,y\rangle\vargeq0\,\}\]
et on pose pour tout $\boldit s\in\interieur C+iM$,
\[L_q(\boldit s,M,C)=
\sum_{y\in C\dual\cap M\dual}q^{-\langle y,\boldit s\rangle}.\]
Par d\'efinition du c\^one dual cette s\'erie converge absolument sur
le c\^one ouvert $\interieur C+iM$ et si $m$ appartient \`a l'int\'erieur de $C$,
la fonction qui \`a $s$ associe $L_q(sm,M,C)$ a un p\^ole d'ordre $\rg M$
en $s=0$.

On pose
\[\alphaetoileV=(\log q)^t\lim_{s\to 1}
(s-1)^tL_q(s\antican,\Pic (\VV),\CeffV)\]
o\`u $t$ d\'esigne le rang de $\Pic (\VV)$
et
\[\alphaV=\frac{\alphaetoileV}{(t-1)!}.\]
Enfin, comme Batyrev et Tschinkel, on pose
\[\betaV=\cardinal H^1(\ff,\Pic(\Vbarre))\]
bien que ce terme soit trivial dans les cas consid\'er\'es ici.
\end{defi}

\begin{rem}
La fonction $L_q(\boldit s,M,C)$ est p\'eriodique de groupe de p\'eriode
contenant $(2\pi i/\log q)M$. Il r\'esulte de la d\'efinition que $\alphaetoileV$
d\'epend du r\'eseau $\Pic(\VV)$, du c\^one $\CeffV$
et de l'\'el\'ement $\antican$ mais
est ind\'ependant de $q$.
\end{rem}

\begin{nota}
La fonction caract\'eristique du c\^one $\CeffV$ est d\'efinie par
\[\forall\boldit s\in\interieur\CeffV,\quad
\chiCeff(\boldit s)=\int_{\CeffV\dual}e^{-\langle s,y\rangle}\Haar y.\]
\end{nota}

\begin{defi}
Si $f$ est une fonction m\'eromorphe sur un ouvert d'un $\CC$-espace vectoriel $W$,
nous dirons que $f$ admet une expression rationnelle en des puissances de $q$
s'il existe une base $(\chi_i)_{1\leq i\leq n}$ de $W\dual$
et une fonction rationnelle
\[R\in\CC(T_1,\dots,T_n)\]
telle que pour tout $\boldit s$ en lequel $f$ est d\'efini on ait
\[f(\boldit s)=R(q^{\langle \chi_1,\boldit s\rangle},
\dots,q^{\langle \chi_n,\boldit s\rangle}).\]
\end{defi}

\begin{prop}
Avec les hypoth\`eses pr\'ec\'edentes, la fonction $L_q(\cdot,\Pic(\VV),
\CeffV)$ admet une expression rationnelle en des puissances de $q$ et
\[\alphaetoileV=\chiCeff(\antican)\in\QQ^*.\]
\end{prop}

\begin{rem}
La raison pour laquelle nous avons substitu\'e $L_q$ \`a $\chiCeff$
dans ce cadre est, qu'\'etant p\'eridique comme la fonction z\^eta
des hauteurs lorsque le syst\`eme de hauteurs v\'erifie la propri\'et\'e (P),
$L_q(\cdot,\Pic(\VV),\CeffV)$ devrait avoir m\^eme lieu singulier que $\zetaHU$
au voisinage de $\antican+i\Pic(\VV)\otimes\RR$.
\end{rem}

\begin{proof}
Rappelons qu'un c\^one $C$ de $M\otimes\RR$ est dit r\'egulier
s'il est de la forme
\[\sum_{i=1}^rm_i\]
o\`u $(m_i)_{1\varleq i\varleq r}$ est une famille libre de $M$.
Pla\c cons-nous tout d'abord dans le cas o\`u
$C$ est une c\^one r\'egulier d'int\'erieur non vide.
La fonction $L_q(\cdot,M,C)$ s'\'ecrit alors
\[L_q(\boldit s,M,C)=\prod_{i=1}^n(1-q^{\langle m_i,\boldit s\rangle})^{-1}\]
et la premi\`ere assertion est imm\'ediate. En outre, si
$\boldit s=\sum_{i=1}^ns_im_i$
alors par \cite[proposition 2.4.5]{batyrevtschinkel:toric},
on a l'\'egelit\'e
\[\chi_C(\boldit s)=\prod_{i=1}^n\frac{1}{s_i}\]
la deuxi\`eme assertion en d\'ecoule aussit\^ot.
Dans le cas g\'en\'eral (cf. \cite[p. 23]{oda:convex}),
on \'ecrit $C\dual$ comme support d'un \'eventail r\'egulier $\Sigma$,
c'est-\`a-dire que $\Sigma$ est un ensemble de c\^ones poly\'edriques 
rationnels strictement convexes de $M\dual\otimes\RR$
tels que:
\begin{itemize}
\item[(i)]si $\sigma\in\Sigma$ et si $\sigma'$ est une face de $\sigma$,
alors $\sigma'\in\Sigma$,
\item[(ii)]si $\sigma,\sigma'\in\Sigma$, alors $\sigma\cap\sigma'$ 
est une face de
$\sigma$ et de $\sigma'$.
\item[(iii)]$C\dual=\cup_{\sigma\in\Sigma}\sigma$,
\item[(iv)]tout $\sigma$ de $\Sigma$ est r\'egulier.
\end{itemize}
Alors si on note $\Sigma^{(i)}$ l'ensemble des \'el\'ements de $\Sigma$
de dimension $i$ et $n=\dim M$, on a les \'egalit\'es:
\[L_q(\boldit s,M,C)=\sum_{\sigma\in\Sigma^{(n)}}L_q(\boldit s,M,\sigma\dual)\]
et
\[\chi_C(\boldit s)=\sum_{\sigma\in\Sigma^{(n)}}\chi_C(\boldit s)\]
et les deux assertions r\'esultent du cas pr\'ec\'edent.
\end{proof}

\subsection{Expression de la constante}

Nous allons maintenant d\'efinir la constante qui appara\^\i t comme
r\'esidu de la fonction z\^eta des hauteurs.

\begin{defi}
On pose
\[\thetaetoileV=\alphaetoileV\betaV\tauV\]
et
\[\thetaV=\alphaV\betaV\tauV\]
\end{defi}

\subsection{G\'eom\'etrie des vari\'et\'es de drapeaux g\'en\'eralis\'ees}

Dans la suite nous nous int\'eressons \`a la situation suivante:

\begin{notas}
Dans la suite $G$ d\'esigne un groupe alg\'ebrique lin\'eaire lisse
semi-simple  et connexe sur $\ff$, $P$ un $\ff$-sous-groupe
parabolique lisse de $G$. On note $\VV$
le quotient $P\sous G$ et $\pi:G\to \VV$
la projection canonique. Par \cite[proposition 2.24]{boreltits:complements},
le rev\^etement universel $\tilde G$ de $G$ est d\'efini sur $\ff$.
Quitte a remplacer $G$ par $\tilde G$ et 
$P$ par son image inverse dans $\tilde G$,
on peut donc supposer que $G$ est simplement connexe.
\par
Pour tout groupe alg\'ebrique lin\'eaire $H$ sur $\ff$,
on note $\Lie(H)$ l'alg\`ebre de Lie restreinte de $H$,
$\rad H$ son radical et $\raduni H$ son radical unipotent.
Le groupe des caract\`eres de $H$ sur $\ff$ est d\'efini par:
\[\caract(H)_\ff=\Hom_{\Spec \ff}(H,\Gmde\ff)\]
et le groupe de cocaract\`eres par
\[\cocaract(H)_\ff=\Hom_{\Spec\ff}(\Gmde\ff,H).\]
\par
On note $P_0$ un $\ff$-sous-groupe parabolique lisse minimal
de $G$ contenu dans $P$. On note $T$ un tore maximal de $\rad 0$
et $S$ sa composante scind\'ee. On a donc les inclusions
\[S\subset T\subset P_0\subset P\subset P\subset G.\]
On fixe \'egalement un sous-groupe de Borel $B$ de
$G^s$ tel que $T^s\subset B\subset P_0^s$.
\par
On note $\racines$ (resp. $\Fracines$, $\positive$, $\Fpositive$)
les racines de $T^s$ (resp. $S$, $T^s$, $S$) dans
$G^s$ (resp. $G$, $B$, $P_0$), $\base$ le base de $\Phi$
associ\'ee \`a $\positive$ et $\Fbase$ celle de $\Fracines$
correspondant \`a $\Fpositive$.
L'application de restriction de $T$ \`a $S$ induit une application
\[j:\base\to\Fbase\cup\{0\}\]
(cf. \cite[\S21.8]{borel:groups}) dont l'image contient
$\Fbase$. le groupe de Weyl de $\racines$ (resp. $\Fracines$)
est note $\Weyl$ (resp. $\FWeyl$). Pour tout $\alpha$ de $\base$
(resp. $\Fbase$), on note $\corac\alpha$ la coracine
correspondante et $\poids\alpha$ le poids fondamental associ\'e.
\par
Pour tout partie $J$ de $\Fbase$, on note $\FWeyl_J$ le sous-groupe
de $\FWeyl$ engendr\'epar les $s_\alpha$, pour $\alpha\in J$
et $\relatif P_J$ le $\ff$ sous-groupe parabolique
\[\relatif P_J=P_0\FWeyl_JP_0.\]
On obtient ainsi une bijection entre les parties de $\Fbase$
et les $\ff$ sous-groupes paraboliques de $G$ contenant $P_0$
avec
\[P_0=\relatif P_\emptyset\qquad\text{et}\qquad G=\relatif P_{\Fbase}.\]
On note $\relatif I$ (resp. I) la partie de $\Fbase$ (resp. $\base$)
correspondant \`a $P$. On a donc
\[I=j^{-1}(\relatif I\cup\{0\})\]
\par
Par \cite[proposition 6.10]{sansuc:brauer}, on a une suite exacte
\[0\to\FF[\VV]\etoile/\ff\etoile\to\ff[G]\etoile/\ff\etoile\to\caract(P)_\ff
\to\Pic(\VV)\to\Pic(G)\]
Il r\'esulte du th\'eor\`eme de Rosenlicht, 
\cite[theorem 3]{rosenlicht:toroidal} que $\ff[G]\etoile/\ff\etoile$
est isomorphe au groupe $\caract(G)_\ff$ et donc trivial
et de \cite[lemma 6.9 (iii)]{sansuc:brauer} que $\Pic(G)$
est nul puisque $G$ est suppos\'e simplement connexe. 
On a donc un isomorphisme canonique
\[\phi:\caract(P)_\ff\iso\Pic(\VV)\]
qui peut \^etre d\'ecrit de la mani\'ere suivante:
si $\chi$ est un caract\`ere de $P$ sur $\ff$, $\chi$ peut \^etre 
vu comme fibr\'e en droites sur $\Spec\ff$ muni d'une action de $P$
et comme fibr\'e en droites sur $\VV$, $\phi(\chi)$ est d\'efini
comme le produit restreint $G\times^P\chi$. Autrement
dit $\phi(\chi)$ est la classe du faisceau $\LL_\chi$
d\'efini par
\[\Gamma(U,L_\chi)=\{\,f\in\Gamma(\pi^{-1}(U),\anneau_G)\mid
\forall y\in\pi^{-1}(U)(\fbarre),\,\forall p\in P(\fbarre),\,
f(pg)=\chi(p)f(g)\}\]
\par
Par \cite[lemme 6.2.10]{peyre:fano}, l'image de $\caract(P)_\fs$
dans $\caract(T)_\fs$ co\"\i ncide avec le sous r\'eseau
engendr\'e par la famille $(\varpi_\alpha)_{\alpha\in\base-I}$
et le c\^one des diviseurs effectifs est donn\'e par
\[\Ceffde\Vs=\sum_{\alpha\in\base-I}\Rplus(-\poids\alpha).\]
L'image de $\caract(P)_\ff$ dans $\caract(T)_\fs$ a donc pour base
la famille
\[\Biggl(\sum_{\beta\in j^{-1}(\alpha)}
\poids\beta\Biggr)_{\alpha\in\Fbase-\relatif I}\]
et $\CeffV$ est donn\'e par
\begin{equation}
\label{equ:cone:effectif}
\sum_{\alpha\in\Fbase-\relatif I}\Rplus\Bigl(
-\sum_{\beta\in j^{-1}(\alpha)}\poids\beta\Bigr).
\end{equation}
\par
Pour tout $J\subset\Fbase$, on note $\relatif\frakr_J$
le radical de $\Lie(\relatif P_J)$ et on pose 
$\frakr=\relatif\frakr_{\relatif I}$. La repr\'esentation adjointe d\'efinit
une action de $P$ sur $\frakr$ et le fibr\'e cotangent $\Omega^1_{\VV/\ff}$
est isomorphe au fibr\'e $G\times^P\frakr$ associ\'e. En prennt
les puissances ext\'erieures maximales, on obtient que
\[\phi(\det\frakr)=\omegaV.\]
On note
\[\demiracines P=\frac{1}{2}\det\frakr\in\caract(P)_\ff\otimes\QQ\]
L'image de $\demiracines P$ dans $\caract (S)_\ff$ par la restriction
co\"\i ncide avec la demi-somme des racines de S compt\'ees avec
des multiplicit\'ees \'egales \`a la dimension de leur espace propre dans
$\frakr$.
\par
Notons en outre que tout choix de bases des sous-espaces propres de 
$\relatif\frakr_\emptyset$ pour l'action de $S$ d\'efinit un isomorphisme
de $\ff$-espace vectoriel $\det\frakr\iso\ff$ et donc
un isomorphisme
\[\omegaV\iso\LL_{2\demiracines P}.\] 
\end{notas}

\begin{rem}
Il r\'esulte des descriptions de $\CeffV$ et de $\antican$
que $\antican\in\interieur\CeffV$ et $\VV$ v\'erifie l'hypoth\`ese (i) 
des hypoth\`eses g\'eom\'etriques ainsi qu l'hypoth\`ese (iv).
L'hypoth\`ese (ii) r\'esulte de \cite[p. 575]{kempf:representations}
et, du fait que $P$ est suppos\'e r\'eduit, la condition (iii)
r\'esulte de la description du groupe du Picard et la derni\`ere 
du fait que $\Vbarre$ est rationnelle. La vari\'et\'e v\'erifie
donc l'ensemble des hypoth\`eses g\'eom\'etriques.
\end{rem}

\subsection{Hauteurs  sur les vari\'et\'es de drapeaux}

Comme dans \cite{fmt:fano}, nous allons tout d'a\-bord
nous restreindre au cas des hauteurs d\'efinies par des sous-groupes
compacts maximaux. Le but de ca paragraphe est d'en rappeler la
construction.

\begin{notas}
Par la d\'ecomposition d'iwasawa (cf. \cite[\S3.3.2]{tits:reductive}),
Il existe pour toute place $\pp$ de $\ff$ un sous-groupe 
compact maximal $\compact_\pp$
de $G(\fp)$ tel que
\begin{equation}
\label{equ:Iwasawa}
G(\fp)=P_0(\fp)\compact_\pp.
\end{equation}
En outre si $\modelede G$ est un mod\`ele de $G$
sur un ouvert de $\CC$, alors par \cite[\S3.9.1]{tits:reductive},
on peut supposer pour presque toute place $\pp$ de $\ff$
que
\[\compact_\pp=\modelede G(\Op)\]
On pose $\compact=\prod_{\pp\in\Mf}\compact_\pp$.
Pour tout caract\`ere $\chi$ de $P$, pour toute place $\pp$
de $\ff$, on consid\`ere la m\'etrique $\metrique$ sur $\LL_\chi$
d\'efinie de la mani\`ere suivante : si $U$ est un voisinage
ouvert de $x$, $\sectiondeL$ une section de $\LL_\chi$
non nulle en $x$ et $\widetilde\sectiondeL$ l'\'el\'ement de
$\Gamma(\pi^{-1}(U),\structural_G)$ qui lui correspond,
\[\forall k\in\compact_\pp,\qquad
\pi(k)=x\Rightarrow\metriquede{\sectiondeL(x)}
=\normede{\widetilde\sectiondeL(k)}.\]
L'existence d'un tel $k$ est assur\'e par \eqref{equ:Iwasawa},
et le terme de droite est in d\'ependant de $k$ puisque pour tout morphisme
continu de $P(\fp)\cap\compact_\pp$ dans $\Rplusetoile$ est trivial.
\par
Le syst\`eme de m\'etriques $\metriques$ ainsi d\'efini est bien ad\'elique.
En effet soit $\modelede P$ l'adh\'erence de $P$ dans $\modelede G$.
Quitte \`a augmenter $\SSS$, on peut supposer que $\modelede P$ est un
sous-groupe parabolique de $\modelede G$ sur $\CC-\SSS$ et le quotient
$\modelede P\sous\modelede G$, bien d\'efini par 
\cite[proposition 1.2]{demazure:parabolique}, fournit un mod\`ele $\mV$
de $\VV$ sur $\CC-\SSS$. On consid\`ere alors le faisceau $\modelede\LL_\chi$
sur $\mV$ d\'efini par
\[\Gamma(\modelede U\!,\modelede\LL_\chi)=
\{\,f\in\Gamma(\pi^{-1}(\modelede U,\structural_\modelede G)\mid
\forall g\in\pi^{-1}(\modelede U)(\fbarre),
\forall p\in P(\fbarre),
f(pg){=}\chi(p)f(g)\,\}\]
pour tout ouvert $\modelede U$ de $\mV$. Quitte \`a augmenter $\SSS$,
on a que $\modelede\LL_\chi$ est un fibr\'e en droites et
un mod\`ele de $\LL_\chi$ et pour tout $\pp\in\Mf-\SSS$
tel que $\compact_\pp=\modelede G(\Op)$ et tout $x$ de $\VV(\fp)$
se relevant en un \'el\'ement $k$ de $\compact_\pp$, la $\Op$-structure
de $\LL_\chi(x)$ d\'efinie par $\LL_\chi$ est induite par la
$\Op$-structure de $\structural_G(c)$ induite par $\structural_\modelede G$,
ce qui mongtre que $\metrique$ co\"\i ncide avec la m\'etrique d\'efinie
par $\modelede\LL_\chi$.
\par
L'application
\[\chi\mapsto(\LL_\chi,\metriques)\]
d\'efinit un syst\`eme de m\'etriques ad\'eliques sur $\VV$
qui v\'erifie de surcro\^\i t la propi\'et\'e (P).
Nous noterons $\HK$
ce syst\`eme de hauteurs. Nous omettrons $\compact$ dans
cette notation lorsqu'aucune confusion ne sera possible.
\end{notas}

\subsection{Enonc\'e du r\'esultat}
Nous pouvons maintenant \'enocer notre r\'esultat proncipal:
\begin{theo}
Avec les notations qui pr\'ec\`edent, la fonction z\^eta
des hauteurs $\zetaHKV(\boldit s)$ converge absolument pour
\[\boldit s\in\antican+\interieur\CeffV+i\Pic\VV\otimes\RR\]
et s\'etend \`a $\Pic\VV\otimes\CC$ en une fonction m\'eromorphe
qui admet une expression rationnelle en des puissance de $q$.
En outre la fonction m\'eromorphe sur $\CC$ qui \`a $s$
associe $\zetaHKV(s\antican)$ a un p\^ole d'ordre $t=\rg\Pic\VV$
en $s=1$ avec
\[\lim_{s\to1}(s-1)^t\zetaHKV(s\antican)=\theta^*(\VV).\]
\end{theo}
\section{D\'emonstration du r\'esultat}
\subsection{Fonction z\^eta et s\'eries d'Eisenstein}

Comme dans \cite[\S2]{fmt:fano}, la d\'emonstration est bas\'ee
sur le fait que la s\'erie z\^eta des hauteurs co\"\i ncide
avec une s\'erie d'Eisenstein, ce qui permet d'appliquer les r\'esultats
d\'emontr\'es par Morris dans \cite{morris:eisenstein:cusp}
et \cite{morris:eisenstein:general} pour ces s\'eries.

\begin{notas}
Pour toute partie $J$ de $\Fbase$, on note $S_J$ le tore
\[\left(\bigcap_{\alpha\in J}\ker\alpha\right)^\circ\]
o\`u pour tout groupe alg\'ebrique $H$, on d\'esigne par $H^\circ$
sa composante neutre. On pose \'egalement
\[\relatif M_J=Z_G(S_J)\qquad\textinmath{et}\quad\relatif N_J=
\raduni(\relatif P_J).\]
Le groupe parabolique $\relatif P_J$
est alors le produit semi-direct de $\relatif M_J$
par $\relatif N_J$. La restriction iduit un isomorphisme
\[\caract(\relatif P_J)\iso\caract(\relatif M_J).\]
Dans la suite, on posera
\[\caractC_J=\caract(\relatif M_J)\otimes_\ZZ\CC
@i\textinmath{Res} ii\caract(S_J)\otimes_\ZZ\CC\]
et on notera $M$ (resp. $N$, $\caractC$, $M_0$, $N_0$,
$\caractCZ$) pour $\relatif M_{\relatif I}$ (resp.
$\relatif N_{\relatif I}$,
$\caractC_{\relatif I}$, 
$\relatif M_\emptyset$, $\relatif N_\emptyset$,
$\caractC_\emptyset$).
\par
Pour tout place $\pp$ de $\ff$, on d\'efinit un morphisme
\[\HPp:M(\Fp)\to\caractC\dual\]
par la relation
\[\forall\chi\in\caract(M),\quad
\forall z\in Z_M(\Fp),\quad\normede{\chi(z)}=
\qq^{\langle\HPp(z),\chi\rangle}\]
et $\HP:M(\adeles)\to\caractC\dual$ est d\'efini comme la somme
des morphisme $\HPp$. On \'etend $\HP$ en une application de
$G(\adeles)$ dans $\caractC\dual$ de la mani\`ere suivante:
\[\forall n\in N(\adeles),\quad
\forall m\in M(\adeles),\quad\forall k\in\compact,\quad
\HP(nmc)=\HP(c).\]
\par
Pour tout sous-groupe compact ouvert $\compact'$ de $\compact$,
on notera $\Cont(P,K')$ l'ensemble des fonctions continues
\[\varphi:P(\adeles)\sous G(\adeles)\to\CC\]
telles que $\varphi$ soit $K'$-finie \`a droite (i.e. l'espace vectoriel
engendr\'e par les translat\'es de $\varphi$ par les \'el\'ements de $K'$
est de dimension finie).
Si $\varphi\in\Cont(P,\compact')$ et $\xi\in\caractC$, on d\'efinit
\[\forall g\in G(\adeles),\quad T_\xi\varphi(g)=\qq^
{\langle \HP(g),\xi\rangle}\varphi(g).\]
\end{notas}

\begin{defi}
Si $\varphi\in\Cont(P,\compact')$ et $\xi\in\caractC$,
la s\'erie d'Eisenstein associ\'ee \`a $\varphi$ et $\xi$
est d\'efinie par
\[\forall g\in G(\adeles),\quad
\EGP(\varphi,\xi,g)=\sum_{\gamma\in P(\ff)\sous G(\ff)}
T_{\xi+\rho_P}\varphi(\gamma g).\]
D'apr\`es un lemme de Godement 
\cite[\S2.2.2, lemme, p. 118]{morris:eisenstein:cusp}
cette s\'erie converge uniform\'ement sur les sous-ensembles
compacts de $G(\adeles)$ d\'es que
$\reel(\xi-\rho_P)$ appartient \`a $C_P$ o\`u $C_P$
est la chambre de Weyl d\'efinie par
\[\forall \alpha\in\Fbase-\relatif I,\quad
(\lambda,\corac\alpha)\vargeq 0.\]
On notera $\EGP(\xi,\cdot)$ la fonction $\EGP(1,\xi,\cdot)$.
\end{defi}

\begin{prop}
\label{prop:zeta:eisenstein}
La s\'erie d\'efinissant la fonction z\^eta des hauteurs $\zetaHV(\boldit s)$
co\"\i ncide avec celle d\'efinissant $\EGP(\boldit s-\rho_P,e)$.
\end{prop}

\begin{proof}
Par \cite[proposition 20.5]{borel:groups}, l'application
$\pi:G(\ff)\to\VV(\ff)$ est surjective. IL suffit donc de v\'erifier
que pour tout \'el\'ement $g$ de $G(\ff)$, on a
\[\HK(\chi)(\pi(g))=\qq^{\langle\HP(g),\chi\rangle}.\]
On \'ecrit donc $g=nmk$ avec $n\in N(\adeles)$, $m\in M(\adeles)$
et $k\in\compact$. Soit $\sectiondeL$ une section de $\LL_\chi$
sur un voisinage ouvert $\UU$ de $\pi(g)$, non nulle
en $\pi(g)$ et correspondant \`a un \'el\'ement $\tilde\sectiondeL$
de $\Gamma(G,\structural_G)$. Par d\'efinition on a
\[\HK(\chi)(\pi(g))=\prod_{\pp\in\Mf}\metriquede{\sectiondeL(\pi(g))}
=\prod_{\pp\in\Mf}\normede{\tilde\sectiondeL(k_\pp)}^{-1}.\]
Mais par la formule du produit 
$\prod_{\pp\in\Mf}\normede{\tilde\sectiondeL(g)}=1$ et donc
\[\Bigl(\prod_{\pp\in\Mf}\normede{\chi(m_\pp n_\pp)}\Bigr)
\prod_{\pp\in\Mf}\normede{\tilde\sectiondeL(k_\pp)}=1.\]
On en d\'eduit les \'egalit\'es
\[\HK(\chi)(\pi(g))=\prod_{\pp\in\Mf}\normede{\chi(m_\pp)}
=\prod_{\pp\in\Mf}\qq^{\langle\HPp(g),\chi\rangle}
=\qq^{\langle\HP(g),\chi\rangle}.\]
\end{proof}

\begin{cor}
La fonction z\^eta des hauteurs $\zetaHKV$ converge absolument dans le
c\^one ouvert
\[\antican+\interieur{\CeffV}+i\Pic\VV\otimes\RR\]
et s'\'etend en une fonction m\'eromorphe sur $\CC$ qui admet
une expression rationnelle en des puissance de $\qq$.
\end{cor}

\begin{proof}
D'apr\`es \cite[\S12.12]{boreltits:reductifs}, on a
\begin{equation}
\label{equ:cone:weyl}
\forall\alpha,\alpha'\in\Fbase,\quad
\alpha \neq\alpha'\Rightarrow\Bigl
\langle\sum_{\beta\in j^{-1}(\alpha)}\poids\beta,\corac\alpha\Bigr\rangle>0
\quad\textinmath{et}\quad
\Bigl\langle\sum_{\beta\in j^{-1}(\alpha)}\poids\beta,
\corac\alpha '\Bigr\rangle=0
\end{equation}
Par cons\'equent le c\^one $C_P$ co\"\i ncide d'apr\`es
\eqref{equ:cone:effectif} avec le c\^one $\CeffV$. La premi\`ere
assertion r\'esulte donc de \cite[lemma, p. 118]{morris:eisenstein:cusp}
mentionn\'e cit\'e ci-dessus. La seconde es un r\'esultat de Morris
\cite[\S6.6, lemma, p. 1164]{morris:eisenstein:general}
\end{proof}

\subsection{Ordre du p\^ole au sommet du c\^one}
L'objectif de ce paragraphe est de montrer le r\'esultat suivant:
\begin{prop}
\label{prop:ordre:pole}
Le lieu des singularit\'es de la fonction $\zetaHKV$
au voisinage du point $s=\antican$ co\"\i ncide avec la r\'eunion
es hyperplans
\[\langle\corac\alpha,\lambda-\antican\rangle=0\]
pour $\alpha\in\Fbase-\relatif I$, chacun de ces hyperplans intervenant
avec une multiplicit\'e au plus \'egale \`a un.
\end{prop}

\begin{rems}
\begin{listrems}
\item Nous verrons plus loin que la multiplicit\'e est en r\'ealit\'e
exactement \'egale \`a un.
\item
Le fait que les singularit\'es soient hyperplanes r\'esulte
de la proposition \ref{prop:zeta:eisenstein}
et de \cite[\S6.6, lemma, p. 1164]{morris:eisenstein:general}.
\end{listrems}
\end{rems}

Comme dans \cite[\S2]{fmt:fano} le principe de la d\'emonstration
est de consid\'erer la fibration
\[P_0\sous P\to P_0\sous G\to \VV,\]
d'exprimer le terme constant des s\'eries d'Eisenstein
correspondant \`a $P_0\sous P$ et $P_0\sous G$ en termes des op\'erateurs
d'entrelacements. Comme dans Harder \cite[p. 278]{harder:chevalley}
ou Morris \cite[\S4.3.4]{morris:eisenstein:cusp}, l'\'etude des singularit\'es
des s\'eries d'Eisenstein se r\'eduit alors \`a la description de celles
des op\'erateurs d'entralecemant qui se d\'eduisent des \'equations
fonctionnelles et du cas de l'op\'erateur associ\'e \`a une r\'eflexion.

\begin{notas}
Pour tout $\varphi$ de $\Cont(P_0,\compact')$ et $\xi$ de $\caractC$,
la s\'erie d'Eisenstein partielle $\EPPZ$ est d\'efinie par
\[\EPPZ(\varphi,\xi,g)=\sum_{\gamma\in P_0(\ff)\sous P(\ff)}
T_{\xi+\rho_{P_0}}\varphi(\gamma g).\]
Pour tout $w\in\FWeyl$ repr\'esent\'e par un \'el\'ement
$w'$ de $\norma_G(S)(k)$, la fonction $C(w,\xi)\varphi$ est
d\'efinie par
\[\forall g\in G(\adeles),\quad C(w,\xi)\varphi(g)=
\int\limits_{\qquad\mathclap{w'N_0(\adeles){w'}^{-1}\cap
N_0(\adeles)\sous N_0(\adeles)}\qquad}
\qq^{\langle\HPZ({w'}^{-1}ng),\xi+\rho_{P_0}\rangle}
\varphi({w'}^{-1}ng)\Haar n\]
o\`u pour tout groupe unipotent $U$ sur $\ff$, la mesure
de Haar sur $U(\adeles)$ est normalis\'ee par
\[\int_{U(\ff)\sous U(\adeles)}\Haar u=1.\]
On note \'egalement $\EPPZ(\xi,g)=\EPPZ(1,\xi,g)$ et
$c(w,\xi)=(C(w,\xi)1)(e).$.
\end{notas}
\begin{rem}
Par \cite[\S2.4.8, theorem, p. 134]{morris:eisenstein:cusp}
et \cite[\S4.3.1, p. 167]{morris:eisenstein:cusp},
on a les \'equations fonctionnelles
\begin{align}
c(w_1,w_2\xi)c(w_2,\xi)&=c(w_1w_2,\xi)
\label{equ:entrelac:fonctionnel}
\\
\noalign{et}
c(w,\lambda)\EGPZ(w\lambda,g)&=\EGPZ(\lambda,g).
\label{equ:entrelac:eisenstein}
\end{align}
En outre, par d\'efinition, pour tout $\xi'$
de $\caractC$ invariant par $w$ on a
\begin{equation}
c(w,\xi+\xi')=c(w,\xi)
\label{equ:entrelac:invariant}
\end{equation}
\end{rem}

\begin{lemme}
Si $w\in\FWeyl$, le lieu des singularit\'es de la fonction $c(w,\cdot)$
au voisinage du point $\lambda=\rho_{P_0}$ est \'egale \`a la r\'eunion
des hyperplans
\[\langle\corac\alpha ,\lambda -\rho_{P_0}\rangle\]
o\`u $\alpha $ d\'ecrit l'ensemble $\{\alpha \in\Fbase\mid w\alpha <0\}$.
La multiplicit\'e de chacun de ces hyperplans est exactement \'egale \`a un.
\end{lemme}

\begin{rem}
Dans le cas o\`u $G$ est un groupe de Chevalley, ce r\'esultat
d\'ecoule imm\'ediatement de ceux de Harder \cite[p. 278]{harder:chevalley}.
\end{rem}

\begin{proof}
Comme dans \cite[p. 429]{fmt:fano}, si $w=s_\alpha$ avec $\alpha\in \Fbase$
il r\'esulte de \eqref{equ:entrelac:invariant} que l'hyperplan
d\'efini par
\[\langle\corac\alpha ,\lambda -\rho _{P_0}\rangle=0\]
contient le lieu singulier de $c(s_\alpha,cdot)$ au voisinage de
$\rho_{P_0}$. Le fait que ce p\^ole soit au plus de multiplicit\'e un
r\'esulte de \cite[\S3.5.2, theorem(i)]{morris:eisenstein:cusp}.
Pour montrer que c'est effectivement un p\^ole, il suffit d'\'ecrire
$c(s_\alpha ,\cdot)$ comme produit de facteurs locaux (cf. aussi
la d\'emonstration du th\'eor\`eme \ref{theo:valeur:constante}).
\par
En g\'en\'eral, on \'ecrit
\[w=s_{\alpha_1}\dots s_{\alpha_q}\]
avec $\alpha_i\in\Fbase$ et $q$ minimal. On pose
\[w_j=s_{\alpha_{j+1}}\dots s_{\alpha_q}.\]
Alors l'\'equation fonctionnelle \eqref{equ:entrelac:fonctionnel}
fournit l'\'egalit\'e
\[c(w,\lambda )=c(s_{\alpha_1},w_1\lambda)\dots c(s_{\alpha_q},\lambda).\]
Mais par \cite[Ch. VI, \S 1.1.6, corollaire 2, p. 158]{bourbaki:lie}, on a
que $w_j^{-1}\alpha _j<0$ et que $\{\alpha \in\Fpositive\mid w\alpha <0\}$
co\"\i ncide avec l'ensemble $\{w_j^{-1},1\varleq j\varleq q\}$.
Il r\'esulte alors de \cite[\S8, sublemma, p. 430]{fmt:fano}
que
\[\langle\alpha_j,w_j\rho_{P_0}\rangle\vargeq
\langle\corac\alpha_j,\rho_{P_0}\rangle\]
avec \'egalit\'e uniquement si $w_j^{-1}\alpha_j\in\Fbase$.
Comme dans \cite{fmt:fano}, on d\'eduit alors du cas qui pr\'ec\`ede
que $c(s_{\alpha _j},w_j\cdot)$ est r\'egulier au voisinage
de $\rho_{P_0}$ sauf si $w_j^{-1}\alpha_j\in\Fbase$ auquel
cas la singularit\'e est contenu dans l'hyperplan
\[\langle w_j^{-1}\alpha_j,\lambda-\rho_{P_0}\rangle=0\]
qui est de multiplicit\'e un.
\end{proof}
\begin{proof}[Fin de la d\'emonstration de la 
proposition \ref{prop:ordre:pole}]
Il r\'esulte de \cite[\S4.3.4]{morris:eisenstein:cusp},
que les singularit\'es de $\EGPZ(\cdot,g)$ sont domin\'ees
par celles de $c(\relatif w_{\Fbase},\cdot)$. le r\'esultat est 
donc d\'emontr\'e dans le cas o\`u $P=P_0$.
Dans le cas g\'en\'eral, pour tout $\xi$ de
\[\rho_{P_0}+C_{P_0}+i\caract(S)\otimes\RR\]
le terme constant de $\EPPZ$ est donn\'e, d'apr\`es
la d\'emonstration de \cite[\S2.3.1, lemma, p. 122]{morris:eisenstein:cusp},
par la formule
\[\int_{N_0(\FF)\sous N_0(\adeles)}\EPPZ(\xi,ng)\Haar n=
\sum_{w!in\relatif W_{\relatif I}}c(w,\xi)q^{\langle\HPZ(g),
w\xi+2\rho_P\rangle}\]
En prenant les syst\`emes d'Eisenstein r\'esiduels successifs
(cf \cite[exemple 3.11.2, p. 1130--1132]{morris:eisenstein:cusp})
on obtient pour tout $\xi$ de $\caractC$ la relation
\[\lim_{\atop{\lambda\to 0}{\lambda\in\caractC^{\perp}}}
\Bigl(\prod_{\alpha\in\relatif I}\langle\corac\alpha,\lambda\rangle\Bigr)
\EPPZ(\lambda+\xi+\rho_{P_0},g)=C_P
q^{\langle\HP(g),\xi+2\rho_P\rangle}\]
o\`u $\caractC^\perp=\{\,\lambda\in\caract(S)\otimes\RR\mid\,
\lambda_{\mid\relatif S_I}=0\}$ et la constante $C_P$
est d\'efinie par
\begin{equation}
\label{equ:def:constante}
C_P=\lim_{\lambda\to\rho_{P_0}}
\Bigl(\prod_{\alpha\in\relatif I}\langle\corac\alpha,\lambda
-\rho_{P_0}\rangle\Bigr)c(\relatif w_{\relatif I},\lambda).
\end{equation}
En sommant sur $P(\ff)\sous G(\ff)$, on obtient
\[\lim_{\atop{\lambda\to 0}{\lambda\in\caractC^{\perp}}}
\Bigl(\prod_{\alpha\in\relatif I}\langle\corac\alpha,\lambda\rangle\Bigr)
\EGPZ(\lambda+\xi+\rho_{P_0},g)=
C_P\EGP(\xi+\rho_P,g)\]
et l'assertion pour $\EGP$ d\'ecoule de celle pour
$\EGPZ$ et du fait que $C_P\neq 0$.
\end{proof}

\begin{rem}
Il d\'ecoule de la d\'emonstration pr\'ec\'edente que
que
\[\lim_{\xi\to 0}\Bigl(\prod_{\alpha\in\base_0-\relatif I}
\langle\corac\alpha,\xi\rangle
\Bigr)\EGP(\xi+\rho_P,g)=C_G/C_P\]
et donc par la proposition \ref{prop:zeta:eisenstein},
\begin{equation}
\label{equ:limite:eisenstein}
\lim_{s\to 1}
(s-1)^{\rg\Pic V}\zetaHKV((s-1)\antican)=\Bigl(
\prod_{\alpha\in\Fbase_{\relatif I}}\langle\corac\alpha,2\rho_{P}\rangle^{-1}
\Bigr){C_G\over C_P}
\end{equation} 
\end{rem}

\subsection{Valeur de la constante}
Par la remarque pr\'ec\'edente, pour cl\^ore la d\'emonstration, il
suffit de comparer $C_G/C_P$ et $\thetaetoileV$ ce qui red\'emontrera
du m\^eme coup que $C_G/C_P$ est une constante non nulle et que
les multiplicit\'es des hyperplans $\langle\corac\alpha,\xi\rangle=0$
sont bien \'egales \`a un.
\begin{theo}
\label{theo:valeur:constante}
On a la relation
\begin{equation}
\label{equ:valeur:constante}
\Bigl(
\prod_{\alpha\in\Fbase_{\relatif I}}\langle\corac\alpha,2\rho_{P}\rangle\Bigr)
{C_G\over C_P}=\thetaetoileV
\end{equation}
\end{theo}
Pour montrer ce r\'esultat, il nous faut d'abord \'ecrire l'op\'erateur
d'entrelacement comme produit d'op\'erateurs locaux que nous allons maintenant d\'efinir.

\begin{notas}
Pour toute place $\pp$ de $\ff$, on note $\relatifp S$ un tore
scind\'e maximal de $G_{\fp}$ tel que
\[\relatif S_{\fp}\subset\relatifp S\subset {P_0}_{\fp},\]
on note $\relatifp P_0$ un sous-groupe parabolique minimal de $G_{\fp}$
tel que
\[\relatifp S\subset\relatifp P_0\subset {P_0}_\fp\]
et $\relatifp N$ son radical unipotent. Quitte \`a modifier le choix
de certains des compacts $\compact_\pp$, on peut supposer que
\[G(\fp)=\relatifp P_0(\fp)\compact_\pp.\]
On note $\Fpracines$ le syst\`eme de $\relatifp S$ dans
$G_{\fp}$ et $\Fpbase$ la base de $\Fpracines$ correspondant
\`a $\relatifp P_0$.
\par
Si $J$ est une partie de $\Fpbase$, od d\'efinit comme pr\'ec\'edemment
le sous groupe parabolique $\relatifp P_J$, l'alg\`ebre de Lie
$\relatifp\frakr_J$, le $\CC$-espace vectoriel $\relatifp\caractC_J$,
l'\'el\'ement $\relatifp w_J$ du groupe de Weyl le caract\`ere
$\rho_{\relatifp P_J}$ et la fonction $H_{\relatifp P_J,\pp}$.
On se donne en outre des bases des sous-espaces propres
de $\relatifp\frakr_{\emptyset}$ pour l'action de $\relatifp S$
de sorte que pour toute partie $J$ de $\Fbase$, l'isomorphisme
de $\ff$-espaces vectoriels
\[\Lambda^{\dim\relatif\frakr_J}(\relatif\frakr_J\otimes\Fp)
\iso\Fp\]
induit par ces bases co\"\i ncide avec celui induit par les
bases choisies sur $\ff$.
\par
Si $U$ est un $\fp$ groupe unipotent de $\relatifp P_{\emptyset}$
ou de son oppos\'e, alors ces bases d\'efinissent un isomorphisme de
vari\'et\'e de $U$ sur $\mathbf A^{\dim U}_\fp$, ce qui permet
de normaliser la mesure de Haar sur $U(\fp)$.
\par
Quitte \`a modifier \`a nouveau certains des $\compact_\pp$,
on peut fixer pour tout $w$ de $\relatif W$ des repr\'esentants
$\tilde w$ appartenant \`a $\norma_G(S)(\ff)\cap\compact$.
\end{notas}

\begin{defi}
Pour tout $\lambda$ de $\relatif\caractC_0$, et tout $w$ de
$\relatif W$, on consid\`ere
\[\relatif c_\pp(w,\lambda)=\quad
\int\limits_{\qquad\mathclap{[\tilde w\relatif N(\fp)\tilde w^{-1}
\cap\relatif N(\fp)]
\sous\relatif N(\fp)}\qquad}\exp(\langle H_{P_0,\pp}(\tilde w^{_1}n),
\lambda +\rho_{P_0}\rangle)\Haar n_\pp\]
\end{defi}

\begin{rem}
Le quotient $C_G/C_P$ se met alors sous la forme
\begin{equation}
\label{equ:valeur:decompose}
{C_G\over C_P}={1\over q^{(g-1)\dim V}}\lim_{\lambda\to0}
\prod_{\alpha\in\Fbase-I}\langle\corac\alpha,\lambda\rangle
\prod_{\pp\in\Mf}{\relatif c_\pp(\relatif w_{\Fbase},
\lambda +\rho _{P_0})\over\relatif c_\pp(\relatif w_I,\lambda +\rho _{P_0})}
\end{equation}
\end{rem}

\begin{defi}
Pour tout $\lambda$ de $\relatifp\caractC_0$ et tout $w$ de $\relatifp W$
se relevant en $w'\in\compact_\pp$
\[c_\pp(w,\lambda)=\quad\int\limits_{
\qquad\mathclap{[w'\relatifp N(\fp){w'}^{-1}
\cap\relatifp N(\fp)]
\sous\relatifp N(\fp)}\qquad}\exp(\langle H_{\relatifp P_0,\pp}({w'}^{-1},
\lambda +\rho _{\relatifp P}\rangle)
\Haar n_\pp\]
\end{defi}

\begin{rem}
Casselman donne dans \cite{casselman:principalseries} une expression
explicite pour $c_\pp(w,\lambda )$ en reliant
$\relatif c_\pp(w,\lambda )$ \`a ce terme, on obtiendra
une expression explicite pour ce dernier.
\end{rem}

\begin{lemme}
Le volume de la vari\'et\'e \`a la place $\pp$ v\'erifie
\[\mesp(V(\fp))={\relatif c_\pp(w_{\Fbase},\rho_{P_0})\over
\relatif c_\pp(w_{\relatif I},\rho_{P_0})}\]
\end{lemme}

\begin{proof}
Ceci r\'esulte imm\'ediatement de la d\'emonstration du lemme 6.2.7
dans \cite{peyre:fano}.
\end{proof}

\begin{lemme}
\label{lemme:biquotient}
Avec les notations pr\'ec\'edentes on a la relation
\[{\relatif c_\pp(\relatif w_{\Fbase},\lambda +\rho _{P_0})
\over\relatif c_\pp(\relatif w_{\relatif I},\lambda +\rho _{P_0})}=
{c_\pp(w_{\Fpbase},\Res\lambda +\rho _{\relatifp P})
\over c_\pp(w_{\relatifp I},\Res\lambda +\rho_{\relatifp P})}\]
o\`u $\relatifp I$ d\'esigne la partie de $\Fpbase$ correspondant
\`a $P_{\fp}$.
\end{lemme}

\begin{proof}
Ce lemme se montre comme le lemme 6.2.8 de \cite{peyre:fano}.
\end{proof}

\begin{nota}
On note
\[\lambda_{P_0}=\sum{\alpha \in j^{-1}(\Fbase)}\poids\alpha\in
\Pic(V_0)\iso\caract(P_0)_\ff\subset\caractC_0.\]
\end{nota}

\begin{lemme}
\label{lemme:convergence}
Le produit
\[\prod_{\pp\in\Mf}L_\pp(s,\Pic\Vbarre)^{-1}
{c_\pp(\relatifp w_{\Fpbase},(s-1)\Res\lambda_{\relatif P}
+\rho_{\relatifp P_0})\over
c_\pp(\relatifp w_{\relatif I},
(s-1)\Res\lambda_{\relatif P}+\rho_{\relatifp P_0})}\]
converge absolument au voisinage de $s=1$.
\end{lemme}
\begin{proof}
Comme dans la d\'emonstration du lemme 6.2.12 de \cite{peyre:fano},
cela r\'esulte de l'expression explicite donn\'ee par Casselmann
de $c_\pp(w,\chi)$.
\end{proof}

\begin{lemme}
La constante $\alphaetoileV$ est donn\'ee par la formule
\[\alphaetoileV={\prod_{\alpha\in\Fbase-\relatif I}
\langle\corac\alpha,\lambda_{\relatif P}\rangle
\over
\prod_{\alpha\in\Fbase-\relatif I}
\langle\corac\alpha,2\rho_P\rangle
}.\]
\end{lemme}
\begin{proof}
Ceci r\'esulte imm\'ediatement de la description
de $\CeffV$ donn\'ee par l'expression \eqref{equ:cone:effectif} et du fait,
d\'eja indiqu\'e en \eqref{equ:cone:weyl} que la matrice de changement
de base passant de $(\corac\alpha)_{\alpha\in\Fbase}$ \`a la base
duale de $(\sum_{\beta\in j^{-1}(\alpha )}\poids\beta)_{\alpha \in\Fbase}$
est diagonale.
\end{proof}
\begin{proof}[D\'emonstration du th\'eor\`eme \ref{theo:valeur:constante}]
Par a formule \eqref{equ:valeur:decompose}, on a
\[{C_G\over C_P}={
\prod\limits_{\alpha \in\Fbase-\relatif I}
\langle\corac\alpha,\lambda_{\relatif P_0}\rangle
\over\qq^{(g-1)\dim V}}
\lim_{s\to 1}(s-1)^{\rg\Pic\VV}\prod_{\pp\in\Mf}
{\relatif c_\pp(w_{\Fbase},(s-1)\lambda _{\relatif P_0}+
\rho_{\relatif P_0})\over
\relatif c_\pp(w_{\relatif I},(s-1)\lambda _{\relatif P_0}+
\rho_{\relatif P_0})}.
\]
Pour tout $\pp$ de $\Mf$, on note $\lambda_\pp(s)=L_\pp(s,\Pic\Vbarre)$
et le quotient $C_G/C_P$ s'\'ecrit
\begin{multline*}
{\prod\limits_{\alpha \in\Fbase-\relatif I}
\langle\corac\alpha,\lambda_{\relatif P_0}\rangle
\over\qq^{(g-1)\dim V}}
\lim_{s\to 1}\Biggl[(s-1)^{\rg\Pic\VV}
L(s,\Pic\Vbarre)\\
\prod_{\pp\in\Mf}
\lambda_\pp(s)^{-1}{\relatif c_\pp(w_{\Fbase},(s-1)\lambda _{\relatif P_0}+
\rho_{\relatif P_0})\over
\relatif c_\pp(w_{\relatif I},(s-1)\lambda _{\relatif P_0}+
\rho_{\relatif P_0})}\Biggr].
\end{multline*}
par les lemmes \ref{lemme:biquotient} et \ref{lemme:convergence},
le produit du bas converge absolument au voisinage de $1$ et le quotient 
$C_G/C_P$ se met sous la forme
\[\begin{split}
&\prod\limits_{\alpha \in\Fbase-\relatif I}
\langle\corac\alpha,\lambda_{\relatif P_0}\rangle
{\lim_{s\to 1}(s-1)^{\rg\Pic\VV}L(s,\Pic\Vbarre)\over
\qq^{(g-1)\dim V}}
\prod_{\pp\in\Mf}
\lambda_\pp(1)^{-1}{\relatif c_\pp(w_{\Fbase},\rho_{\relatif P_0})\over
\relatif c_\pp(w_{\relatif I},\rho_{\relatif P_0})}\\
&=\prod\limits_{\alpha \in\Fbase-\relatif I}
\langle\corac\alpha,\lambda_{\relatif P_0}\rangle\tau_{\HK}(\VV)
\end{split}\]
ce qui conclut la d\'emonstration.
\end{proof}

\par
\begin{merci}
Je tiens \`a remercier Laure Blasco pour ses pr\'ecieuses
indications.
\end{merci}

\ifx\undefined\bysame
\newcommand{\bysame}{\leavevmode\hbox to3em{\hrulefill}\,}
\fi
\ifx\undefined\numero
\newcommand{\numero}{$\hbox{n}^\circ$}
\fi
\ifx\undefined\andname
\newcommand{\andname}{and }
\fi
\ifx\undefined\comma
\newcommand{\comma}{,}
\fi

\end{document}